\documentclass[aop,MSNbibl,citesort,seceqn,dvips]{arximspdf}
\usepackage[mathscr]{eucal}
\usepackage{accents}

\doi{10.1214/11-AOP741} 
\volume{41}
\issue{2}
\pubyear{2013}
\firstpage{989}
\lastpage{1029}

\makeatletter
\renewcommand{\mathring}[1]{\accentset{\circ}{#1}}
\newcommand{\eqref}[1]{(\ref{#1})}
\newtheorem{theorem}{Theorem}[section]
\newtheorem{proposition}[theorem]{Proposition}
\newtheorem{lemma}[theorem]{Lemma}

\newproclaim{definition}[theorem]{Definition}
\newproclaim{remark}[theorem]{Remark}
\newproclaim{conjecture}[theorem]{Conjecture}
\newproclaim{notation}{Notation}

\def\E{\mathbb{E}}
\def\N{\mathbb{N}}
\def\P{\mathbb{P}}
\def\Q{\mathbb{Q}}
\def\R{\mathbb{R}}
\def\Z{\mathbb{Z}}

\def\1{\mathbf{1}}

\newcommand{\ltwocv}{\stackrel{L^2}{\longrightarrow}}
\newcommand{\pcv}{\stackrel{P}{\longrightarrow}}
\newcommand{\ascv}{\stackrel{\mathrm{a.s.}}{\longrightarrow}}
\newcommand{\ascvi}{\stackrel{a.s.}{\longrightarrow}}

\def\B{\mathscr{B}}
\def\C{\mathscr{C}}
\def\D{\mathscr{D}}
\def\sE{\mathscr{E}}
\def\F{\mathscr{F}}
\def\G{\mathscr{G}}
\def\M{\mathscr{M}}
\def\Pscr{\mathscr{P}}
\def\T{\mathscr{T}}

\def\O1{\hat{\Omega}}
\def\Fi{\hat{\F} }
\def\Ft{\hat{\F}_t}
\def\Fth{\hat{\F}^*_t}
\def\PT{\hat{\P}_T}
\def\tK{\tilde{K}}
\def\brK{\mathring{K}}
\def\bY{\bar{Y}}
\def\bZ{\bar{Z}}
\def\Lipp{\mathrm{Lip}_1}

\def\Xb{ \tilde{X} }
\def\Zb{\bar{Z}}
\def\Yb{\bar{Y}}
\def\Kb{\tK}
\def\phib{\bar{\phi}}
\def\Tx{\eta}
\def\Xp{X'}
\def\XbD{\Xb^D}
\def\exp{\operatorname{exp}}

\def\supp{\operatorname{Supp}}

\makeatother

\begin{document}
\begin{frontmatter}

\title{A super Ornstein--Uhlenbeck process interacting with its center
of mass}
\runtitle{SOU interacting with its COM}

\begin{aug}
\author[A]{\fnms{Hardeep} \snm{Gill}\corref{}\ead[label=e1]{gillhs@interchange.ubc.ca}}
\runauthor{H. Gill}
\affiliation{University of British Columbia}
\address[A]{Department of Mathematics\\
University of British Columbia\\
1984, Mathematics Road \\
Vancouver, British Columbia\\
V6T1Z2, Canada\\
\printead{e1}} 
\end{aug}

\received{\smonth{8} \syear{2011}}
\revised{\smonth{12} \syear{2011}}

%
\begin{abstract}


We construct a supercritical interacting measure-valued diffusion with
representative particles that are attracted to, or repelled from, the
center of mass. Using the historical stochastic calculus of Perkins, we
modify a super Ornstein--Uhlenbeck process with attraction to its
origin, and prove continuum analogues of results of
Engl\"ander [\textit{Electron. J. Probab.} \textbf{15} (2010)
1938--1970] for binary branching Brownian motion.

It is shown, on the survival set, that in the attractive case the mass
normalized interacting measure-valued process converges almost surely
to the stationary distribution of the Ornstein--Uhlenbeck process,
centered at the limiting value of its center of mass. In the same
setting, it is proven that the normalized super Ornstein--Uhlenbeck
process converges a.s. to a Gaussian random variable, which
strengthens a theorem of Engl\"ander and
Winter [\textit{Ann. Inst. Henri Poincar\'e Probab. Stat.}
\textbf{42} (2006)
171--185] in this
particular case. In the repelling setting, we show that the center of
mass converges a.s., provided the repulsion is not too strong and then
give a conjecture. This contrasts with the center of mass of an
ordinary super Ornstein--Uhlenbeck process with repulsion, which is
shown to diverge a.s.

A version of a result of Tribe
[\textit{Ann. Probab.} \textbf{20} (1992) 286--311] is proven on the
extinction set; that is, as it approaches the extinction time, the
normalized process in both the attractive and repelling cases converges
to a random point a.s.
\end{abstract}

%


\end{frontmatter}


\section{Introduction and main results}\label{sec1}

The existence and uniqueness of a self-interacting measure-valued
diffusion that is either attracted to or repelled from its centre of
mass is shown below. It is natural to consider a super
Ornstein--Uhlenbeck (SOU) process with attractor (repeller) given by the
centre of mass of the process as it is the simplest diffusion of this sort.

This type of model first appeared in a recent paper of Engl\"
ander \cite
{Eng2010} where a $d$-dimensional binary Brownian motion, with each
parent giving birth to exactly two offspring and branching occurring at
integral times, is used to construct a binary branching
Ornstein--Uhlenbeck process where each particle is attracted (repelled)
by the center of mass (COM). This is done by solving the appropriate
SDE along each branch of the particle system and then stitching these
solutions together.

This model can be generalized such that the underlying process is a
branching Brownian motion (BBM), $\T$ (i.e., with a general offspring
distribution). We might then solve an SDE on each branch of $\T$:
%
\begin{equation}\label{Ebps}\quad
Y_i^{n}(t) = Y_{p(i)}^{n-1}(n-1) + \gamma\int_{n-1}^t \bar{Y}_n(s) -
Y_i^{n}(s) \,ds +\int_{n-1}^t dB_i^n(s)\vspace*{-2pt}
\end{equation}
for $n-1< t\le n$, where $B_i^n$ labels the $i$th particle of $\T$
alive from time $n-1$ to $n$, $p(i)$ is the parent of $i$ and
\[
\bar{Y}_n(s) = \frac{1}{\tau_n} \sum_{i=1}^{\tau_n} Y_i^n(s)\vspace*{-2pt}
\]
is the center of mass. Here, $\tau_n$ is the population of particles
alive from time $n-1$ to~$n$. This constructs a branching OU system
with attraction to the COM when $\gamma>0$ and repulsion when $\gamma<0$.

It seems reasonable then, to take a scaling limit of branching particle
systems of this form and expect it to converge in distribution to a
measure-valued process where the representative particles behave like
an OU process attracting to (repelling from) the COM of the process.
Though viable, this approach will be avoided in lieu of a second method
utilizing the historical stochastic calculus of Perkins \cite
{Perkins2002} which is
more convenient for both constructing the SOU interacting with its COM
and for proving various properties. The idea is to use a supercritical
historical Brownian motion to construct the interactive SOU process by
solving a certain stochastic equation. This approach for constructing
interacting measure-valued diffusions was pioneered in \cite
{Perkins1995} and utilized in, for example,~\cite{HG2009}.

A supercritical historical Brownian motion, $K$, is a stochastic
process taking values in the space of measures over the space of paths
in $\R^d$. One can think of~$K$ as a supercritical superprocess which
has a path-valued Brownian motion as the underlying process. That is,
if $B_t$ is a $d$-dimensional Brownian motion, then $\hat{B}_t =
B_{\cdot\wedge t}$ is the underlying process of $K$. More information
about $K$ is provided in Section~\ref{sec2}.

It can be shown that if a path $y\dvtx [0,\infty)\rightarrow\R^d$ is chosen
according to $K_t$ (loosely speaking---this is made rigorous below in
Definition~\ref{Dcamp}), then $y(s)$ is a Brownian motion stopped at
$t$. Projecting down gives
\[
X_t^K(\cdot) = \int\1(y_t \in\cdot)K_t(dy),\vspace*{-2pt}
\]
which is a (supercritical) super Brownian motion.\vadjust{\goodbreak}

A key advantage to projecting down and constructing measure-valued
processes is that it is possible to use the historical stochastic
calculus to couple different projections together (and hence couple
measure-valued diffusions).

One can sensibly define the Ornstein--Uhlenbeck SDE driven
by $y$ according to~$K_t$ as the solution of a stochastic equation.

\begin{definition} Let $Z_0\dvtx \R^d\rightarrow\R^d$ be Borel measurable.
We say that $(X, Z)$ is a solution to the strong equation~\ref{se1} if the pair satisfies
{\renewcommand{\theequation}{$\mbox{(SE)}_{Z_0, K}^1$}
\begin{eqnarray}\label{se1}
\qquad\qquad &&\mbox{(a)}\quad\hspace*{3.5pt}   Z_t( y) = Z_0(y_0)
+y_t-y_0-\gamma\int_0^t Z_s( y) \,ds,\qquad K\mbox{-a.e.},
\nonumber
\\[-10pt]
\\[-10pt]
\nonumber&&
\mbox{(b)}\quad   X_t(A) = \int1(Z_t\in A)K_t(dy)\qquad \forall A\in\B(\R
^d)\ \forall t\ge0,
\end{eqnarray}}
\hspace*{-2pt}where $X$ and $Z$ are appropriately adapted. We will henceforth call
the projection~$X$ the ordinary super Ornstein--Uhlenbeck process.
\end{definition}

If $\gamma>0$, $Z_t$ is attracted to the origin, and if $\gamma<0$ it
is repelled. The approximate meaning of $K$-a.e. in (a) is that the
statement holds for~$K_t$-a.a. $y$, for all~$K_t$, $\P$-a.s. The exact
definition is given in the next section. The projection $X_t$ is a SOU
process with attraction to (repulsion from) the origin at rate $\gamma
$. Intuitively, $K$ tracks the underlying branching structure and $Z_t$
is a function transforming a typical Brownian path into a typical
Ornstein--Uhlenbeck path. 

Note that in the definitions of~\ref{se1} and~\ref{SE2} given below, part (b) is unnecessary to solve the
equations. It has been included to provide an easy comparison to the
strong equation of Chapter V.1 of~\cite{Perkins2002}.

For all the results mentioned in the remainder of this work, the
standing assumption (unless indicated otherwise) will be that
%
\setcounter{equation}{1}
\begin{equation}\label{Efinmass}
\int1 \,dK_0 < \infty.
\end{equation}

\begin{theorem}\label{TSOUrep}
There is a pathwise unique solution $(X, Z)$ to~\ref{se1}. That is, $X$ is unique $\P$-a.s. and $Z$ $K$-a.e. unique.
Furthermore, the map $t\rightarrow X_t$ is continuous and $X$ is a
$\beta$-super-critical super Ornstein--Uhlenbeck process.
\end{theorem}

Similar to the above, we establish a function of the path $y$ that is a
path of an OU process with attraction (repulsion) to the COM and
project down.

\begin{definition}
Let $Y_0\dvtx \R^d\rightarrow\R^d$ be Borel measurable. Define ($X', Y$) as
the solution of
{\renewcommand{\theequation}{$\mbox{(SE)}_{Y_0, K}^2$}
\begin{eqnarray} \label{SE2}
\qquad\qquad&&\textup{(a)}\quad\hspace*{4pt}   Y_t( y) = Y_0(y_0) +y_t
-y_0+\gamma\int_0^t \bar{Y}_s- Y_s( y) \,ds,\qquad K\mbox{-a.e. }
\nonumber
\\[-10pt]
\\[-10pt]
\nonumber&&
\textup{(b)} \quad  X'_t(A) = \int1(Y_t\in A)K_t(dy)\qquad \forall A\in\B(\R
^d)\ \forall t\ge0,
\end{eqnarray}}
\hspace*{-2pt}where the COM is
\[
\bY_s = \frac{\int x X'_s(dx)}{\int1 X'_s(dx)}.
\]
We will call the projection $X'$ the super Ornstein--Uhlenbeck process
with attraction (repulsion) to its COM, or the interacting super
Ornstein--Uhlenbeck process.
\end{definition}

Note that by our definitions of $X$ and $X'$ as solutions to $\mbox{(SE)}_{Y_0, K}^i, i=1,2$, respectively, that \eqref{Efinmass} is the
same as saying that
$ \int1 \,dX_0, \int1 \,dX'_0 <\infty, $
as these quantities equal that in \eqref{Efinmass}.

\begin{theorem}\label{TExUniq}
There is a pathwise unique solution to~\ref{SE2}.
\end{theorem}

One could prove this theorem using a combination of the proof of
Theorem~\ref{TSOUrep} and a localization argument. We find it more
profitable however, to employ a correspondence with the ordinary SOU
process $X$. This correspondence plays a central role in the analysis
of $X'$, and indeed reveals a very interesting structure: We have that
for any $\gamma$,
%
\setcounter{equation}{2}
\begin{equation} \label{Ecoupl}
\int\phi(x) \,dX'_t(x) = \int\phi\biggl(x+ \gamma\int_0^t \Zb_s \,ds\biggr) \,dX_t(x),
\end{equation}
where $\Zb$ is the COM of $X$, defined as $ \Zb_s = \frac{\int x
X_s(dx)}{\int1 X_s(dx)}.$
The correspondence essentially says that the SOU process with
attraction (repulsion) to its COM is the same as the ordinary SOU
process being dynamically pushed by its COM. From this equation, a
relation between $\Yb$ and $\Zb$ can be established:
\[
\Yb_t = \Zb_t + \gamma\int_0^t \Zb_s \,ds.
\]
Define $\Xb_t = \frac{X_t}{X_t(1)}$ and $\Xb'_t = \frac
{X_t'}{X'_t(1)}$. As the goal of this work is to prove that $\Xb'_t$
has interesting limiting behavior as $t$ approaches infinity, \eqref
{Ecoupl} yields a method of approach: show first that the time
integral in \eqref{Ecoupl} converges in some sense and establish
limiting behavior for $\Xb$. One then hopes to combine these two facts
with~\eqref{Ecoupl} to get the desired result. 

Let $S$ be the event that $K_t$ survives indefinitely. Note that this
implies that on $S$, both $X$ and $X'$ survive indefinitely by their
definitions as solutions of the equations above. Let $\Tx$ be the time
at which $K_t$ goes extinct.


The next two theorems settle the question of what happens on the
extinction set~$S^c$.
%
\begin{theorem}\label{TExtinc}
On $S^c$, $\Yb_t$ and $\Zb_t$ converge as $t\uparrow\eta<\infty$,
$\P
$-a.s., for any $\gamma\in\R$.
\end{theorem}

\begin{theorem}\label{TExtConv}
On the extinction set, $S^c$,
\[
\Xb_t\rightarrow\delta_F\quad \mbox{and}\quad \Xb'_t\rightarrow\delta_{F'}
\]
as $t\uparrow\Tx<\infty$ a.s., where $F$ and $F'$ are $\R^d$-valued
random variables such that
\[
F' = F + \gamma\int_0^{\eta} \Zb_s\,ds.
\]
\end{theorem}

This last theorem is an analogue of the result of Tribe \cite
{Tribe1992} for ordinary critical superprocesses. Note that here it
does not matter whether there is attraction or repulsion from the COM.

The following three theorems for the attractive case ($\gamma>0$) form
the main results of this work.

\begin{theorem}\label{TSurvival}
On $S$ the following hold:
\begin{longlist}
\item[(a)] If $\gamma> 0$
\[
\Zb_t \ascvi 0 \quad\mbox{and} \quad\Yb_t \ascvi\gamma\int_0^\infty\Zb_s \,ds,
\]
and this integral is finite almost surely.
\item[(b)] If $\gamma=0$, then $\Zb_t = \Yb_t$ and this quantity
converges almost surely.
\end{longlist}
\end{theorem}

This says that the COMs of ordinary and interacting SOU process
converge and is the result that allows us to fruitfully use the
correspondence of \eqref{Ecoupl} to show convergence of the
interacting SOU process.

The next theorem shows that the mass normalized SOU process converges
almost surely, which is a new result among superprocesses. Engl\"ander
and Winter in~\cite{EW2006} have shown that this process converges in
probability, and before them, Engl\"ander and Turaev in~\cite{ET2002}
shown convergence in distribution.

One expects a result of this sort to hold since in the particle picture
conditional on survival, at large time horizons, there are a very large
number of particles that move as independent OU processes, each of
which are located in the vicinity of the origin. Thus, we expect that
in the limit the mass will be distributed according to the limiting
distribution of an OU process.

\begin{theorem} \label{TSOUconv}
Suppose $\gamma>0$. Then on $S$,
\[
d(\Xb_t, P_\infty)\ascvi0,
\]
where $P_t$ is the semigroup of an Ornstein--Uhlenbeck process with
attraction to $0$ at rate $\gamma,$ and $d$ is the Vasserstein metric
on the space of finite measures on $\R^d$ (see Definition~\ref{DVas}).
\end{theorem}

\begin{remark}
It is possible to show that Theorem~\ref{TSOUconv} holds for a more
general class of superprocesses. If the underlying process has an
exponential rate of convergence to a stationary distribution, then the
above theorem goes through. One can appeal to, for example, Theorem 4.2
of Tweedie and Roberts~\cite{TR2000} for a class of such continuous
time processes.
\end{remark}


Using the correspondence of \eqref{Ecoupl}, one can then show that the
mass normalized interacting SOU process with attraction converges to a
Gaussian distribution, centered at the limiting value of the COM, $\Yb
_\infty$.

\begin{theorem} \label{TISOUConv}
Suppose $\gamma>0$. Then on $S$,
%
\begin{equation} \label{EIMconv}
d(\Xb'_t,P_\infty^{\Yb_\infty}) \ascvi0,
\end{equation}
where $P_\infty^{\Yb_\infty}$ is the OU-semigroup at infinity, with the
origin shifted to $\Yb_\infty$.
\end{theorem}

When there is repulsion, matters become more difficult on the survival
set. It is no longer clear whether there exists a limiting random
measure, or what the correct normalizing factor is. We can show however
that in some cases the COM of the interacting SOU process still
converges. That there should be a limiting measure comes from the fact
that the ordinary repelling SOU process has been shown to converge in
probability by Engl\"ander and Winter in~\cite{EW2006} to a multiple
of Lebesgue measure. One may expect something similar to hold for the
interacting SOU process, given \eqref{Ecoupl}. Unfortunately, the
correspondence is rendered ineffectual in this case by part (b) of the
following theorem.

\begin{theorem}\label{TrCOM}
The following hold on $S$ if $0>\gamma> -\frac{\beta}{2}$:
\begin{longlist}[(a)]
\item[(a)] The process, $\Yb_t$ converges almost surely.
\item[(b)] However, $\Zb_t$ diverges exponentially fast. That is, there
is a random variable~$L$, such that
\[
\P(e^{\gamma t}\Zb_t \rightarrow L, L \neq0 | S) =1.
\]
\end{longlist}
\end{theorem}

This reveals an interesting byplay between $X$ and $X'$ in the
repelling case. That is, if one fixes a compact set $A\subset\R^d$,
then for the ordinary SOU process, $X_t$, $A$ is exponentially distant
from the COM of the process. However, the COM of $X'$ will possibly lie
in the vicinity of $A$ for all time. Therefore, one might expect that
$A$ is charged by a different amount of mass by $X'_t$ than $X_t$, and
thus we might need to renormalize $X'_t$ differently to get a valid
limit. It is also possible that the limit for each case is different
(and not simply connected by a random translation).

The proofs for these theorems and more are contained in the following
sections. In Section~\ref{sec2}, we give some background information for the
historical process $K$ and some rigorous definitions. In Section~\ref{Sexpre}, we
prove Theorems~\ref{TSOUrep} and~\ref{TExUniq} and state some
important preliminary results regarding the nature of the support of a
supercritical historical Brownian motion. These are consequently used
to get moment bounds on the center of mass processes $\Yb_t$ and $\Zb
_t$. We also derive a martingale problem for $\Xb$ and $\Xb'$. In
Section~\ref{SConvergence}, we give the proofs of the convergence theorems mentioned
above and in the final section give the proofs of technical results
that are crucial to prove these. 



\section{Definitions and background material}\label{sec2}

\begin{notation*} We collect some terms below:\vspace*{9pt}

\begin{tabular}{ l l }
$E, E'$& Metric spaces \\
$C_c(E, E')$& Compact, cont. functions from $E$ to $E'$ \\
$C_b(E, E')$& Bounded, cont. functions from $E$ to $E'$ \\
$(E,\sE)$& Arbitrary measure space over $E$\\
$b\sE$& Bounded, $\sE$-mble. real-valued functions\\
$M_F(E)$& Space of finite measures on $E$\\
$\mu(f)= \int f\,d\mu$& where $f\dvtx E\rightarrow\R$ and $\mu$ a measure on
$E$\\
$\mu(f) = (\mu(f_1),\ldots,\mu(f_n))$& if $f= (f_1,\ldots, f_n)$,
$f_i\dvtx E\rightarrow\R$\\
$|p|$& $p\in\R^d$, denotes Euclidean norm of $p$\\
$\| f\| = \sup_{E} \sum\vert f_i(x)\vert$& if $f= (f_1,\ldots, f_n)$,
$f_i\dvtx E\rightarrow\R$\\
$C = C(\R_+, \R^d)$& Space of continuous paths in $\R^d$ \\
$\C$& The Borel $\sigma$-field of $C$ \\
$y^t = y_{\cdot\wedge t}$& The path $y$ stopped at $t$ \\
$C^t = \{y^t\dvtx  y \in C\}$& Set of all paths stopped at $t$\\
$\C_t = \sigma(y^s, s\le t, y \in C)$& Natural filtration associated
with $\C$
\end{tabular}\vspace*{6pt}
\end{notation*}


We take $K$ to be a supercritical historical Brownian motion.
Specifically, let $K$ be a $(\Delta/2, \beta, 1)$-historical
superprocess (here $\Delta$ is the $d$-dimensional Laplacian), where
$\beta>0$ constant, on the probability space $(\Omega, \F, (\F
_t)_{t\ge0}, \P)$. Here $\beta$ corresponds to the branching bias in
the offspring distribution, and the $1$ to the variance of the
offspring distribution. A martingale problem characterizing $K$ is
given below. For a more thorough explanation of historical Brownian
motion than found here, see Section V.2 of~\cite{Perkins2002}.

It turns out that $K_t$ is supported on $C^t\subset C$ a.s. and
typically, $K_t$ puts mass on those paths that are ``Brownian'' (until
time $t$ and fixed thereafter). As $K$ takes values in $M_F(C)$,
$K_t(\cdot)$ will denote integration over the $y$ variable.

Let $\hat{B}_t = B(\cdot\wedge t)$ be the path-valued process
associated with $B,$ taking values in $C^t$. Then for $\phi\in b\C$,
if $s\le t$ let $P_{s, t}\phi(y) = \E^{s,y}(\phi(\hat{B}_t))$, where
the right-hand side denotes expectation at time $t$ given that until
time $s$, $\hat{B}$ follows the path~$y$.

The weak generator, $\hat{A}$, of $\hat{B}$ is as follows. If $\phi\dvtx
\R
_+\times C \rightarrow\R$ we say
$\phi\in\D(\hat{A})$ if and only if $\phi$ is bounded, continuous and
($\C_t$)-predictable, and for some $\hat{A}_s\phi(y)$ with the same
properties as $\phi,$
\[
\phi(t, \hat{B}) - \phi(s, \hat{B}) - \int_s^t \hat{A}_r \phi
(\hat
{B})\,dr,\qquad t\ge s,
\]
is a ($\C_t$)-martingale under $P_{s, y}$ for all $s\ge0, y\in C^s.$

If $m\in M_F(\R^d)$, we will say $K$ satisfies the historical
martingale problem,~\ref{hmp}, if and only if $K_0 = m$ a.s. and
$\forall\phi\in\D(\hat{A}),$
{\renewcommand{\theequation}{$(\mathrm{HMP})_m$}
\begin{eqnarray}\label{hmp}
 M_t(\phi) \equiv K_t(\phi_t) - K_0(\phi_0) - \int_0^t K_s(\hat
{A}_s\phi)\,ds-\beta\int_0^t K_s(\phi_s)\,ds
\nonumber\\\\[-18pt]
\eqntext{\displaystyle\mbox{is a continuous $(\F_t)$-martingale}}\\
 \eqntext{\displaystyle\mbox{with }
\langle M(\phi)\rangle_t = \int_0^t K_s(\phi_s^2)
\,ds\ \forall t\ge0, \mbox{ a.s.}}
\end{eqnarray}}

Using the martingale problem~\ref{hmp}, one can construct an orthogonal
martingale measure $M_t(\cdot)$ with the method of Walsh~\cite{W1984}.
Denote by $\Pscr$, the $\sigma$-field of $(\F_t)$-predictable sets in
$\R_+\times\Omega$. If $\psi\dvtx \R_+\times\Omega\times C
\rightarrow\R$
is $\Pscr\times\C$-measurable and
%
\setcounter{equation}{0}
\begin{equation} \label{Emmeas}
\int_0^t K_s(\psi_s^2) \,ds < \infty\qquad \forall t\ge0,
\end{equation}
then there exists a continuous local martingale $M_t(\psi)$ with
quadratic variation $\langle M(\psi) \rangle_t = \int_0^t K_s(\psi^2_s)
\,ds$. If the expectation of the term in \eqref{Emmeas} is finite, then
$M_t(\psi)$ is an $L^2$-martingale.

\begin{definition}\label{Dcamp}
Let $(\O1,\Fi, \Ft)= (\Omega\times C, \F\times\C, \F_t \times
\C
_t)$. Let $\Fth$ denote the universal completion of $\Ft$. If $T$ is a
bounded $(\F_t)$-stopping time, then the normalized Campbell measure
associated with $T$ is the measure $\PT$ on $(\O1, \Fi)$ given by
\[
\PT(A\times B) = \frac{\P(\1_A K_T(B))}{m_T(1) }\qquad\mbox{for } A\in
\F,
B\in\C,
\]
where $m_T(1) = \P(K_T(1))$. We denote sample points in $\O1$ by
$(\omega, y)$. Therefore, under $\PT$, $\omega$ has law $K_T(1)\,d\P
\cdot m_T^{-1}(1)$ and conditional on $\omega$, $y$ has law $K_T(\cdot
)/K_T(1)$.
\end{definition}

\begin{definition} \label{DKae}
For two $(\Fth)$-measurable processes $Z^1$ and $Z^2,$ we will say that
$Z^1 = Z^2, K$-a.e. if
\[
Z^1(s, \omega, y) = Z^2(s, \omega, y)\qquad \forall s\le t, K_t\mbox
{-a.a. } y
\]
for all fixed times $t\ge0$.
\end{definition}

\begin{definition}
We say that $(X, Z)$ is a solution to the strong equation
\ref{se1} if both (a) and (b) of that equation are satisfied
and where $Z_t$ is an $(\Fth)$-predictable process and $X_t$ is an
$(\F
_t)$-predictable process.

We say $(X', Y)$ is the solution to the stochastic equation
\ref{SE2} if it satisfies (a) and (b) of that equation and
where $Y_t$ is an $(\Fth)$-predictable process, $X'_t$ is an $(\F
_t)$-predictable process.

Unless stated otherwise, $(X, Z)$ will refer to a solution of
\ref{se1} and $(X', Y)$ to the solution of~\ref{SE2}.
\end{definition}

\begin{definition} For an arbitrary $\Fth$-adapted, $\R^d$-valued
process $z_t$, define the centre of mass (COM), $\bar{z}_t$, with
respect to $K_t$ as follows:
\[
\bar{z}_t = \frac{K_t(z_t)}{K_t(1)}.
\]
We also let
\begin{eqnarray*}
\Xb_t(\cdot) &\equiv&\frac{X_t(\cdot)}{X_t(1)} = \frac{K_t(Z_t\in
\cdot
)}{K_t(1)},\qquad
\Xb_t'(\cdot) \equiv\frac{X'_t(\cdot)}{X'_t(1)} =
\frac
{K_t(Y_t\in\cdot)}{K_t(1)}\quad \mbox{and}\\
\tK_t(\cdot) &=& \frac
{K_t(\cdot
)}{K_t(1)}.
\end{eqnarray*}
\end{definition}

Note that as $K$ is supercritical, $K_t$ survives indefinitely on a set
of positive probability $S$, and goes extinct on the set of positive
probability $S^c.$ Hence, we can only make sense of ${\bar{z}}_t$ for
$t<\Tx$
where $\Tx$ is the extinction time. 

\begin{definition}
Let the process $M_t$ be defined for $t\in[0,\zeta)$ where $\zeta\le
\infty$ possibly random. $M$ is called a local martingale on its
lifetime if there exist stopping times $T_N \uparrow\zeta$ such that
$M_{T_N\wedge\cdot}$ is a martingale for all $N$. The interval $[0,
\zeta)$ is called the lifetime of $M$.
\end{definition}

The following definition introduces a metric on the space of finite
measures, which is equivalent to the metric of convergence in
distribution on the space of probability measures.

\begin{definition}\label{DVas}
Let $d$ denote the Vasserstein metric on the space of finite measures
on $\R^d$. That is, for $\mu, \nu$ finite measures,
\[
d(\mu, \nu) = \sup_{\phi\in\Lipp}\int\phi(x) \,d(\mu- \nu)(x),
\]
where $\Lipp= \{\psi\in C(\R^d)\dvtx  \forall x, y, \vert\psi(x) -
\psi (y)\vert\le\vert x-y\vert, \|\psi\| \le1\}$.
\end{definition}



\section{Proofs of existence and preliminary results}\label{Sexpre}
\mbox{}
\begin{pf*}{Proof of Theorem~\ref{TSOUrep}}
Although closely related, this does not follow automatically from
Theorem 4.10 of Perkins~\cite{Perkins1995} where it is shown that
equations like $\mbox{(SE)}^2$ [but with more general, interactive,
drift and diffusion terms in (a)] have solutions if $K$ is a \textit
{critical} historical Brownian motion.

Note that $\brK_t = e^{-\beta t}K_t$ defines a $(\Delta/2, 0,
e^{-\beta
t})$-Historical superprocess. Let $\PT^1$ be the Campbell measure
associated with $\brK$ (note that if $T$ is taken to be a random
stopping time, this measure differs from $\PT$). The proof of Theorem~V.4.1 of~\cite{Perkins2002} with minor modifications shows that
$\mbox{(SE)}_{Z_0, \brK}^1$\vspace*{2pt} has a pathwise unique solution. This is because (K3) of Theorem
2.6 of~\cite{Perkins1995} shows that under $\PT^1$, $y_t$ is a
Brownian motion stopped at time $T$ and Proposition~2.7 of the same
memoir can be used to replace Proposition~2.4 and Remark~2.5(c)
of~\cite{Perkins2002}
 for the setting where the branching variance depends on time.

Once this is established, it is simple to deduce that if $({\mathring
{X}}, Z)$ is
the solution of $\mbox{(SE)}_{Z_0, \brK}^1$, and we let
\[
X_t(\cdot) \equiv e^{\beta t}{\mathring{X}}_t(\cdot) = \int1(Z_t\in
\cdot)
K_t(dy)
\]
then $(X,Z)$ is the pathwise unique solution of~\ref{se1}. The only thing to check this is that $Z_t(\omega, y) =
Z_0(\omega
, y_0) +y_t-y_0-\gamma\int_0^t Z_s(\omega, y) \,ds $ $K$-a.e., but this
follows from the fact that $\brK_t \ll K_t, \forall t$.

It can be shown by using by Theorem 2.14 of~\cite{Perkins1995} that
$X$ satisfies the following martingale problem: For $\phi\in C^2_b(\R^d)$,
\[
M_t(\phi) \equiv X_t(\phi) - X_0(\phi) + \int_0^t \int\gamma
x\cdot
\nabla\phi(x) -\frac{\Delta}{2}\phi(x) X_s(dx)\,ds - \beta\int
_0^tX_s(\phi
) \,ds
\]
is a martingale where $\langle M(\phi)\rangle_t = \int_0^tX_s(\phi
^2)\,ds$. Then by Theorem II.5.1 of~\cite{Perkins2002} this implies
that $X$ is a
version of a SOU process, with initial distribution given by
$K_0(Z_0^{-1}(\cdot))$.
\end{pf*}

\begin{remark} \label{RSOUHrep}
(a) Under the Lipschitz assumptions of Section V.1 of
\cite{Perkins2002}, one can in fact uniquely solve
\ref{se1} where
$K$ is a supercritical historical Brownian motion. The proof above can
be extended with minor modifications.\vspace*{-6pt}
\begin{longlist}[(b)]
\item[(b)] The proof of Theorem~\ref{TSOUrep} essentially shows that
under $\PT$, $T$ fixed, the path process $y\dvtx \R^+\times\O1
\rightarrow
\R^d$ such that $(t, (\omega, y)) \mapsto y_t$ is a $d$-dimensional
Brownian motion stopped at $T$.
\item[(c)] Under $\PT$, $T$ fixed, $Z_t$ can be written explicitly as a
function of the driving path: For $t\le T$,
%
\begin{eqnarray}\label{RZSDE}
e^{\gamma t}Z_t &=& Z_0 + \int_0^t e^{\gamma s} \,dZ_s + \int_0^t Z_s
(\gamma e^{\gamma s}) \,ds
\nonumber
\\[-8pt]
\\[-8pt]
\nonumber
&=&Z_0 +\int_0^t e^{\gamma s} \,dy_s + \int_0^t e^{\gamma s} (-\gamma Z_s)
\,ds + \int_0^t Z_s (\gamma e^{\gamma s}) \,ds,
\end{eqnarray}
where we have used a differential form of~\ref{se1}(a)
for the second equality. Hence,
%
\begin{equation}\label{EZpath}
Z_t(y) = e^{-\gamma t}Z_0 + \int_0^t e^{-\gamma(t-s)}\,dy_s.
\end{equation}
\end{longlist}
\end{remark}

Next, we show that there exists a unique solution to~\ref{SE2}. 

\begin{pf*}{Proof of Theorem~\ref{TExUniq}}
Suppose there exists a solution $Y$ satisfying~\ref{SE2}. Then under $\PT$, $Y_t$ can be written as a function of the
driving path $y$ and $\bY$. Using integration by parts gives
\begin{eqnarray*}
e^{\gamma t} Y_t &=& Y_0 + \int_0^t e^{\gamma s} \,dY_s + \int_0^t
\gamma
e^{\gamma s} Y_s \,ds \\
&=& Y_0 + \int_0^t e^{\gamma s} \,dy_s + \int_0^t \gamma e^{\gamma s}
\bar
{Y}_s \,ds
\end{eqnarray*}
and hence
%
\begin{equation}\label{RYSDE}
Y_t =e^{-\gamma t}Y_0 + \int_0^t e^{\gamma(s-t)} \,dy_s + \int_0^t
\gamma e^{\gamma(s-t)} \bar{Y}_s \,ds.
\end{equation}
If $(X, Z)$ is the solution to~\ref{se1} where $Z_0 =
Y_0$, then note that by Remark~\ref{RSOUHrep}(c),
\[
Y_t =Z_t+ \int_0^t \gamma e^{\gamma(s-t)} \bar{Y}_s \,ds.
\]

By taking the normalized measure $\tK_t$ on both sides of the above
equation, we get
%
\begin{equation} \label{Ecorr}
\bY_t = \bar{Z}_t + \gamma\int_0^t e^{-\gamma(t-s)} \bar{Y}_s \,ds.
\end{equation}
Hence, $\bY_t$ is seen to satisfy a Volterra Integral Equation of the
second kind (see equation (2.2.1) of~\cite{PM2008}) and therefore can be
solved pathwise to give
\[
\bY_t = \bZ_t + \gamma\int_0^t \bZ_s \,ds,
\]
which is easily verified using integration by parts. Also, if $\bY^1_t$
is a second process which solves \eqref{Ecorr}, then
\begin{eqnarray*}
|\bY_t -\bY^1_t| &=& \biggl|\gamma\int_0^t e^{-\gamma(t-s)} (\bar
{Y}_s -
\bY^1_s) \,ds\biggr| \\
&\le&|\gamma|\int_0^t e^{-\gamma(t-s)} |\bar{Y}_s - \bY^1_s| \,ds.
\end{eqnarray*}
By Gronwall's inequality, this implies $\bY_t = \bY^1_t$, for all $t$
and $\omega$. Pathwise uniqueness of $\Xp$ follows from the uniqueness
of the solution to~\ref{se1} and the uniqueness of the
process $\bY_t$ solving \eqref{Ecorr}.

We have shown that if there exists a solution to~\ref{SE2}, then it is necessarily pathwise unique. Turning now to existence
to complete the proof, we work in the opposite order and define $Y$ and
$\Xp$ as functions of the pathwise unique solution to
\ref{se1} where $Z_0 = Y_0$:
\begin{eqnarray*}
Y_t& =& Z_t + \gamma\int_0^t \bZ_s \,ds, \\
X'_t(\cdot) &=& K_t(Y_t\in\cdot).
\end{eqnarray*}
Then $\Yb_t$ satisfies the integral equation \eqref{Ecorr}, and hence
\[
\int_0^t \Zb_s \,ds = \int_0^t e^{-\gamma(t-s)}\Yb_s \,ds.
\]
Therefore
\begin{eqnarray*}
Y_t &=& Z_t + \gamma\int_0^t e^{-\gamma(t-s)}\Yb_s \,ds \\
&=& e^{-\gamma t}Y_0 + \int_0^t e^{-\gamma(t-s)} \,dy_s + \gamma\int_0^t
e^{-\gamma(t-s)}\Yb_s \,ds,
\end{eqnarray*}
by equation (\ref{EZpath}), and so
\[
e^{\gamma t} Y_t = Y_0 + \int_0^t e^{\gamma s} \,dy_s + \gamma\int_0^t
e^{\gamma s}\Yb_s \,ds.
\]
Multiplying by $e^{-\gamma t}$ and using integration by parts shows
\begin{eqnarray*}
Y_t = Y_0 + y_t-y_0 + \gamma\int_0^t (\Yb_s - Y_s)\,ds
\end{eqnarray*}
which holds for $K$-a.e. $y$, thereby showing $(X',Y)$ satisfies
\ref{SE2}.
\end{pf*}

\begin{remark} \label{RCOMrel}
Some useful equivalences in the above proof are collected below. If
$Y_0 = Z_0$, then for $t<\Tx$,
\begin{eqnarray*}
&&\mbox{(a)}\hspace*{68pt}\quad Y_t = Z_t + \gamma\int_0^t\Zb_s\,ds, \\
&&\mbox{(b)}\hspace*{67pt}\quad \Yb_t = \Zb_t + \gamma\int_0^t\Zb_s\,ds, \\
&&\mbox{(c)}\hspace*{45pt}\quad  Y_t - \bY_t  = Z_t - \bZ_t, \\
&&\mbox{(d)}\quad \int_0^t e^{-\gamma(t-s)}\bY_s \,ds = \int_0^t \bZ_s\,ds.
\end{eqnarray*}
\end{remark}

These equations intimately tie the behaviour of the interacting and
ordinary SOU processes. Part (a) says that the interacting SOU process
with attraction to the center of mass is the same as the ordinary SOU
process pushed by the position of its center of mass.

We now consider the martingale problem for $X'$. For $\phi\dvtx  \R
^d\rightarrow\R$, recall that $\bar{\phi}_t \equiv K_t(\phi
(Y_t))/K_t(1)$ and that the lifetime of the process $\bar{\phi}$ is
$[0, \Tx)$. Then the following theorem holds:

\begin{theorem}\label{TIto}
For $\phi\in C^2_b(\R^d, \R),$ and $t< \eta$,
\[
\phib_t = \phib_0 + N_t+ \int_0^t \bar{b}_s \,ds,
\]
where
\[
b_s = \gamma\nabla\phi(Y_s)\cdot(\Yb_s - Y_s) + \tfrac{1}{2}\Delta
\phi(Y_s)
\]
and $N_t$ is a continuous local martingale on its lifetime such that
\[
N_t = \int_0^t\int\frac{\phi(Y_s) - \phib_s}{K_s(1)}\,dM(s, y)
\]
and hence has quadratic variation given by
\[
[N]_t = \int_0^t\frac{\overline{\phi_s^2} - (\phib_s)^2}{K_s(1)}\,ds.
\]
\end{theorem}

Similarly, the following is true.

\begin{remark} \label{RIto2}
The method of Theorem~\ref{TIto} can be used to show that for the
$\beta$-supercritical SOU process, X, for $\phi\in C^2_b(\R^d, \R)$
and $t<\eta$,
\[
\Xb_t(\phi) = \Kb_t(\phi(Z_t)) = \Kb_0(\phi(Z_0)) + N_t+ \int
_0^t \Kb
_s(L\phi(Z_s)) \,ds,
\]
where
\[
L\phi(x) = -\gamma x\cdot\nabla\phi(x) + \tfrac{1}{2}\Delta\phi(x)
\]
and $N_t$ is a continuous local martingale on its lifetime such that
\[
N_t = \int_0^t\int\frac{\phi(Z_s) - \Kb_s(\phi
(Z_s))}{K_s(1)}\,dM(s, y)
\]
and hence has quadratic variation given by
\[
[N]_t = \int_0^t\frac{\Kb_s(\phi^2(Z_s)) - \Kb_s(\phi
(Z_s))^2}{K_s(1)}\,ds.
\]
\end{remark}

\begin{pf*}{Proof of Theorem~\ref{TIto}}
The proof is not very difficult; one need only use It\^{o}'s lemma
followed by some slight modifications of theorems in Chapter V of~\cite
{Perkins2002} to deal with the drift introduced in the historical martingale problem
due to the supercritical branching.

Let $T$ be a fixed time and $t\le T$. Recall that under $\PT$, $y$ is a
stopped Brownian motion by Remark~\ref{RSOUHrep}(b), and hence
$Y_t(y)$ is a stopped OU process (attracting to~$\Yb_t$). Therefore,
under $\PT$:
\begin{eqnarray*}
\phi(Y_t) &=& \phi(Y_0) + \int_0^t \nabla\phi(Y_s)\cdot dY_s+ \frac
{1}{2}\sum_{i,j\le d} \int_0^t \phi_{ij}(Y_s)\,d[Y^i, Y^j]_s \\
&=&\phi(Y_0) + \int_0^t \nabla\phi(Y_s)\cdot dy_s+ \int_0^t \gamma
\nabla
\phi(Y_s)\cdot(\Yb_s-Y_s)+ \frac{1}{2}\Delta\phi(Y_s)\,ds
\end{eqnarray*}
by the classical It\^{o}'s lemma. Then
\begin{eqnarray*}
K_t(\phi(Y_t)) &=& K_t(\phi(Y_0)) + K_t\biggl(\int_0^t \nabla\phi
(Y_s)\cdot dy_s\biggr)+K_t\biggl(\int_0^t b_s \,ds\biggr)\\
&=& K_0(\phi(Y_0))+ \int_0^t\int\phi(Y_0)\,dM(s, y)+ \beta\int_0^t
K_s(\phi(Y_0))\,ds \\
&&{}+ \int_0^t \int\biggl[ \int_0^s \nabla\phi(Y_r)\cdot dy_r
\biggr]\,dM(s,y)\\
&&{}+\beta\int_0^t K_s\biggl[ \int_0^s \nabla\phi(Y_r)\cdot
dy_r\biggr]\,ds\\
&&{}+\int_0^t\int b_s\, dM(s,y)+ \beta\int_0^t K_s(b_s)\,ds +\int
_0^tK_s(b_s)\,ds\\
&=& K_0(\phi(Y_0))+ \int_0^t\int\phi(Y_s)\,dM(s, y)\\
&&{}+ \beta\int_0^t K_s(\phi(Y_s))\,ds + \int_0^tK_s(b_s)\,ds.
\end{eqnarray*}
The equality in the third line to $K_t(\int_0^t \nabla\phi
(Y_s)\cdot dy_s)$ follows from Proposition~2.13 of \cite
{Perkins1995}. The equality of the fourth line to the last term in the
first line follows from a generalization of Proposition V.2.4 (b)
of~\cite{Perkins2002}. The last equality then follows by collecting like terms and using the
definition of $\Yb$.

Note that $K_t(1)$ is Feller's $\beta$-supercritical branching
diffusion and hence
\[
K_t(1) = K_0(1) + M^0_t + \beta\int_0^t K_s(1) \,ds,
\]
where $M^0_t$ is a martingale such that $[M^0]_t = \int_0^t K_s(1) \,ds.$

Therefore for $t<\Tx$, It\^{o}'s formula and properties of $K_t(1)$ and
$ K_t(\phi(Y_t))$ imply that
\begin{eqnarray*}
\phib_t &=& \frac{K_t(\phi(Y_t))}{K_t(1)} \\
&=& \phib_0 + \int_0^t\int\biggl[\frac{\phi(Y_s)}{K_s(1)}- \frac
{K_s(\phi
(Y_s))}{K_s(1)^2} \biggr]\,dM(s,y)+\int_0^t\frac{K_s(b_s)}{K_s(1)}\,ds.
\end{eqnarray*}
Since $\phi$ is bounded, the stochastic integral term can be localized
using the stopping times $T_N \equiv\min\{t\dvtx  K_t(1)\ge N \mbox{ or }
K_t(1)\le1/N\}\wedge N$ and hence it is a local martingale on $[0,\Tx
)$. It is easy to check that it has the appropriate quadratic variation.
\end{pf*}

The following lemmas will be used extensively in Section \ref
{SConvergence}, but will be proven in Section~\ref{Stech}.

\begin{lemma} \label{LExtProb}
There is a nonnegative random variable $W$ such that
\[
e^{-\beta t} K_t(1) \rightarrow W\qquad \mbox{a.s.}
\]
and $\{\Tx< \infty\} = \{W=0\}$ almost surely.
\end{lemma}

Note that as $X$ and $X'$ are defined as projections of $K$, their mass
processes are the same as that of $K$, and thus grow at the same rate.

\begin{definition} Let
\[
h(\delta) = \bigl(\delta\ln^+(1/\delta)\bigr)^{1/2},
\]
where $\ln^+(x) =\ln x \vee1$. Let
\[
S(\delta, c) = \{y\dvtx  |y_t - y_s| < ch(|t-s|), \forall t,s \mbox{ with }
|t-s|\le\delta\}.
\]
\end{definition}

\begin{lemma}\label{LCSP}
Let $K$ be a supercritical historical Brownian motion, with drift~$\beta$, branching variance $1$, and initial measure $X_0$. For $c_0>6$
fixed, $c(t) = \sqrt{t + c_0}$, there exists a.s. $\delta(\omega)>0$
such that $\supp(K_t(\omega))\subset S(\delta(\omega), c(t))$ for all~$t$. Further, given $c_0$, $\P(\delta<\lambda) < p_{c_0}(\lambda)$
where $p_{c_0}(\lambda)\downarrow0$ as $\lambda\downarrow0$ and for
any $\alpha>0$, $c_0$ can be chosen large enough so that
$p_{c_0}(\lambda) = C(d,c_0)\lambda^\alpha$ for $\lambda\in[0,1]$.
\end{lemma}


The following moment estimates are useful in establishing the
convergence of~$\bar{Y}_t$. Recall that $\Tx$ is the extinction time of
$K$. 

\begin{lemma}\label{Lsecmom}
Assume $\P(\tK_0(|Y_0|^2+|y_0|^2))<\infty$. Then,
\[
\P(\overline{\vert Y_{t}\vert^2}; t< \Tx) < A(\gamma, t),
\]
where
\begin{eqnarray*}
A(\gamma, t) =
\cases{
O(1+t^6e^{-2\gamma t}),&\quad $\mbox{ if } \gamma< 0$, \vspace*{2pt}\cr
 O(1+t^5),&\quad $\mbox{ if } \gamma\ge0.$
}
\end{eqnarray*}
\end{lemma}

\begin{remark} \label{Rmoments}
(a) The proof of Lemma~\ref{Lsecmom}, under the same hypotheses
(if $Z_0=Y_0$) yields
\[
\P(\overline{\vert Z_{t}\vert^2};t<\Tx) < A(\gamma, t).\vspace*{-6pt}
\]
\begin{longlist}[(b)]
\item[(b)] Lemma~\ref{Lsecmom} and its proof can be extended to show
that for any positive integer $k$ if $\P(\tK_0(|Z_0|^k+
|y_0|^k))<\infty$ (and $Z_0 = Y_0$), then there exists a function
$B(\gamma, t, k)$ polynomial in $t$ if $\gamma\ge0$, exponential if
$\gamma<0$ such that
\[
\P(\overline{\vert Y_{t}\vert^k}; t<\Tx) < B(\gamma, t, k) \quad\mbox{and}\quad \P
(\overline{\vert Z_{t}\vert^k}; t<\Tx)< B(\gamma, t, k).
\]
\end{longlist}
\end{remark}


\section{Proofs of convergence} \label{SConvergence}

We will henceforth, unless specified otherwise, assume that for a path
$y \in C$, $Z_0(y_0) = Y_0(y_0) = y_0$. Recall that, by the
construction of solutions to $\mbox{(SE)}^i,$
\[
\int\phi(y_0)K_0(dy) = \int\phi(x) X_0(dx) = \int\phi(x) X'_0(dx).
\]
Also recall that our standing hypothesis is that $K_0$ has finite
initial mass [and hence $X_0'(1)=X_0(1)<\infty$].

In this section, we will first settle what happens to $X$ and $X'$ when
there is extinction, and then the case when there is survival, under
the attractive regime and lastly address the interacting repelling SOU
process on the survival set.


\subsection{On the extinction set}
\mbox{}
\begin{pf*}{Proof of Theorem~\ref{TExtinc}}
Assume for now that $\P(\tK_0(\vert y_0\vert^2))<\infty$ (the case
where $K_0 = 0$ can be ignored without loss of generality). By
Theorem~\ref{TIto},
\[
\Yb_t = \Yb_0 + \int_0^t\int\frac{Y_s - \Yb_s}{K_s(1)}\,dM(s, y)
\]
and therefore is a local martingale on its lifetime with reducing
sequence $\{T_N\}$ as defined in the proof of the same theorem.\vadjust{\goodbreak} Using
Doob's weak inequality and Lemma~\ref{Lsecmom},
\begin{eqnarray*}
\P\Bigl(\sup_{s< t\wedge\Tx} |\Yb_s| > n\Bigr) &=&\lim_{N
\rightarrow
\infty} \P\Bigl(\sup_{s< t\wedge T_N} |\Yb_s| > n\Bigr) \\
&\le&\lim_{N\rightarrow\infty} \frac{1}{n^2}\E(|\Yb_{t\wedge T_N}|^2)
\\
&<& \frac{A(\gamma, t)}{n^2}.
\end{eqnarray*}
By the first Borel--Cantelli lemma,
\[
\P\Bigl(\sup_{s< t\wedge\Tx} |\Yb_t| > n \mbox{ i.o.}\Bigr) = 0.
\]
It follows that
\[
\liminf_{s\rightarrow t\wedge\Tx} \Yb_s >-\infty\quad\mbox{and}\quad
\limsup
_{s\rightarrow t\wedge\Tx} \Yb_s < \infty
\]
which implies that on the set $\{\Tx< t\}$, $\Yb_s$ converges, by
Theorem IV.34.12 of Rogers and Williams~\cite{RW1985}. This shows
convergence on the extinction set as $S^c = \bigcup_t \{\Tx< t\}$.

Note that if $\nu(\cdot) = \P(K_0 \in\cdot)$. Theorem II.8.3
of~\cite{Perkins2002}
 gives
\[
\P(K\in\cdot) = \int\P_{K_0}(K\in\cdot)\,d\nu(K_0),
\]
where $\P_{K_0}$ is the law of a historical Brownian motion with
initial distribution $\delta_{K_0}$. Hence, the a.s. convergence of
$\Yb_t$ in the case where $\tK_0(|y_0|^2)$ is finite in mean imply
%
\begin{eqnarray}\label{Egen1}
\P\Bigl(\lim_{t\uparrow\Tx} \Yb_t \mbox{ exists} ; S^c\Bigr)&= &\int\P
_{K_0}\Bigl(\lim_{t\uparrow\Tx} \Yb_t \mbox{ exists} ; S^c \Bigr)\,d\nu
(K_0) \nonumber\\
&=& \int\P_{K_0}(S^c)\,d\nu(K_0)\\
&=&\P(S^c) \nonumber
\end{eqnarray}
if $K_0(|y_0|^2+1)<\infty, \nu$-a.s.

Finally, to get rid of the assumption that $K_0(|y_0|^2)<\infty$
note that Corollary 3.4 of~\cite{Perkins1995} ensures that if $K_0(1)
< \infty$, then at any time $t>0$, $K_t$ (and hence $X_t, X'_t$) is
compactly supported. Therefore, letting $S_r = \{K_r \neq0\}$ we
see that
\begin{eqnarray*}
\P\Bigl(\lim_{t\uparrow\Tx}\Yb_t \mbox{ does not exist}, S^c\Bigr) &=&
\P\biggl(\bigcup_{r\in\N}\Bigl\{\lim_{t\uparrow\Tx}\Yb_t \mbox{ does not
exist}, S_{1/r}, S^c\Bigr\}\biggr)\\
&\le&\sum_{r\in\N} \P\Bigl(\P_{K_{1/r}} \Bigl(\lim_{t\uparrow\Tx }\Yb_t
\mbox{ does not exist}, S^c \Bigr)\1_{S_{1/r}} \Bigr) \\
&=& 0
\end{eqnarray*}
by (\ref{Egen1}) since $K_{1/r}$ a.s. compact implies that
$K_{1/r}(|y_{1/r}|^2)<\infty$ holds. This completes the proof for the
convergence of $\Yb$ on $S^c$ in its full generality.

The convergence of $\Zb_t$ now follows from the convergence of $\Yb_t$
and equation~(\ref{Ecorr}).
\end{pf*}


\begin{pf*}{Proof of Theorem~\ref{TExtConv}}
As in the previous proof, note that we need only consider the case that
$\P(\tK_0(\vert y_0\vert^2))<\infty$.

We will follow the proof of Theorem 1 of Tribe~\cite{Tribe1992} here. Define
\[
\zeta(t) = \int_0^t\frac{1}{K_s(1)}\,ds,\qquad t<\Tx.
\]
It is known by the work of Konno and Shiga~\cite{KS1988} in the case
where $\beta=0$, that $\zeta\dvtx [0,\Tx)\rightarrow[0,\infty)$
homeomorphically (recall that $\Tx<\infty$ a.s. in that case). This
latter result also holds when $\beta>0$ on the extinction set $S^c$ by
a Girsanov argument.

Define $D\dvtx [0,\zeta(\Tx-))\rightarrow[0,\Tx)$ as the unique inverse of
$\zeta$ (on $S^c$, this defines the inverse on $[0,\infty)$) and for
$t\ge\zeta(\Tx-)$, let $D_t = \infty$. Let
\[
X^D_t = X'_{D_t},\qquad \Xb^D_t = \frac{X^D_t}{X^D_t(1)}\quad \mbox{and}\quad\G_t
= \F_{D_t}
\]
and define
\[
L_t\phi(x) = \gamma(\Yb_t - x)\cdot\nabla\phi(x)+\tfrac
{1}{2}\Delta\phi(x).
\]
Let
\[
T_N=\int_0^{\eta_N} \frac{1}{K_s(1)}\,ds,
\]
where $\eta_N = \inf\{s\dvtx  K_s(1) \le1/N\}$. Then note that $T_N
\uparrow\zeta(\eta-)$ and each $T_N$ is a $\G_t$-stopping time.
On $S^c$ for $\phi\in C^2_b$, Theorem~\ref{TIto} implies
%
\begin{eqnarray} \label{Eexteq}
\XbD_{t\wedge T_N}(\phi) &=& \Xb_0(\phi) + \int_0^{D_{t\wedge T_N}}
\Xb
'_s(L_s\phi)\,ds + M_{D_{t\wedge T_N}}(\phi)\nonumber\\
&=& \Xb_0(\phi) + \int_0^{t\wedge T_N} \XbD_s(L_{D_s}\phi
)X^D_s(1)\,ds +
N_{t\wedge T_N}(\phi)\\
&= &\Xb_0(\phi) + \int_0^{t\wedge T_N} X^D_s(L_{D_s}\phi)\,ds +
N_{t\wedge
T_N}(\phi)\nonumber
\end{eqnarray}
since $dD_t = X^D_t(1)\,dt$ and where $N_t = M_{D_t}$. It follows that
$N_{t\wedge T_N}$ is a $\G_t$-local martingale. Then, by Theorem~\ref{TIto},
\begin{eqnarray*}
[\XbD(\phi)]_{t\wedge T_N} &=& \int_0^{D_{t\wedge T_N}} \frac
{\Xb
'_s(\phi^2) - \Xb'_s(\phi)^2}{X_s(1)}\,ds \\
&=& \int_0^{t\wedge T_N} \XbD_s(\phi^2) - \XbD_s(\phi)^2 \,ds,
\end{eqnarray*}
which is uniformly bounded in $N$. Hence, sending $N\rightarrow\infty
$, one sees that $N_{t\wedge\zeta(\Tx-)}$ is a $\G_t$-martingale.

Note that on $S^c$, $\zeta(\Tx-) = \infty$ and hence on that event,
\begin{eqnarray*}
\int_0^\infty X^D_s(|L_{D_s} \phi|) \,ds &\le&\int_0^\infty\int\vert
\gamma(\bY_{D_s}-x)\cdot\nabla\phi(x)\vert+\frac{1}{2}\vert
\Delta \phi (x)\vert X^D_s(dx)\,ds \\
&\le&\int_0^\infty|\gamma|\|\nabla\phi\|K_{D_s}(|\bY
_{D_s}-Y_{D_s}|) +\frac{1}{2}\|\Delta\phi\|X^D_s(1)\,ds \\
&=& \int_0^\infty X^D_s(1)\biggl( |\gamma|\|\nabla\phi\|\tK
_{D_s}(|\bY _{D_s}-Y_{D_s}|) +\frac{1}{2}\|\Delta\phi\|\biggr)\,ds \\
&\le&\int_0^\infty X^D_s(1)\biggl( |\gamma|\|\nabla\phi\|
(\overline {|Y_{D_s}|^2})^{{1}/{2}} +\frac{1}{2}\|\Delta\phi\|
\biggr)\,ds,
\end{eqnarray*}
where in the second line we have used the definition of $X'$ and the
Cauchy--Schwarz inequality in the fourth. Using the definition of $D_s$ yields
%
\begin{eqnarray}\label{Eextbd}
\int_0^\infty X^D_s(|L_{D_s} \phi|) \,ds &\le&\biggl( |\gamma|\|\nabla
\phi
\|\sup_{s<\Tx}(\overline{|Y_{s}|^2})^{{1}/{2}} +\frac
{1}{2}\|
\Delta\phi\| \biggr) \int_0^\infty X^D_s(1)\,ds \nonumber\\
&=& \biggl( |\gamma|\|\nabla\phi\|\sup_{s<\Tx}(\overline
{|Y_{s}|^2})^{{1}/{2}} +\frac{1}{2}\|\Delta\phi\| \biggr) \Tx
\\
&<&\infty\nonumber
\end{eqnarray}
as $\phi\in C_b^2$ and $\overline{\vert Y_s\vert^2}$ is continuous
on $[0,
\Tx)$ (which follows from Theorem~\ref{TIto}). Hence, this implies
that for $\phi$ positive, on $S^c$
\[
N_t(\phi) > -\Xb_0(\phi) - \int_0^\infty X^D_s(|L_{D_s} \phi|) \,ds
\]
for all $t$ and hence by Corollary IV.34.13 of~\cite{RW1985}, $N_t$
converges as $t\rightarrow\infty$. Therefore by \eqref{Eexteq} and
\eqref{Eextbd}, $\Xb^D_t(\phi)$ converges a.s. as well.

Denote by $\XbD_\infty(\phi)$ the limit of $\XbD_t(\phi)$. It is
immediately evident that $\XbD_\infty(\cdot)$ is a probability measure
on $\R^d$. To show that $\XbD_\infty(\cdot) = \delta_{F'}$ where $F'$
is a random point in $\R^d$, we now defer to the proof of Theorem 1 in
Tribe~\cite{Tribe1992}, as it is identical from this point forward.

Similar (but simpler) reasoning holds to show $\Xb_t \rightarrow
\delta
_F$ a.s. on $S^c$ where $F$ is a random point in $\R^d$. Let $f(t) =
\gamma\int_0^t \Zb_s \,ds$. Note that $f$ is independent of $y$ and that
$f(t) \rightarrow f(\eta)$ a.s. when $t\uparrow\eta$ because $\Yb_t =
\Zb_t + f(t)$ and both $\Yb_t$ and $\Zb_t$ converge a.s. by
Theorem \ref
{TExtinc}. Then for $\phi$ bounded and Lipschitz,
\begin{eqnarray*}
\biggl\vert\int\phi\bigl(x - f(t)\bigr) \Xb'_t(dx) - \int\phi\bigl(x - f(\eta)\bigr) \Xb
'_t(dx) \biggr\vert & \le& C\vert f(\eta) - f(t)\vert \\
&\ascv&0
\end{eqnarray*}
as $t\uparrow\eta$. Therefore it is enough to note that since $f(\eta)$
depends only on $\omega$, the convergence of $\Xb'_t$ gives
\[
\int\phi\bigl(x - f(\eta)\bigr) \Xb'_t(dx) \ascv\phi\bigl(F' - f(\eta)\bigr)
\]
and hence
\[
\int\phi\bigl(x - f(t)\bigr) \Xb'_t(dx) \ascv\phi\bigl(F' - f(\eta)\bigr).
\]
By Remark~\ref{RCOMrel}(a),
\begin{eqnarray*}
\int\phi\bigl(x - f(t)\bigr) \Xb'_t(dx) &=& \int\phi\bigl(Y_t(y)- f(t) \bigr)\tK_t(dy)
\\
&=& \int\phi(Z_t(y))\tK_t(dy) \\
&=& \int\phi(x)\Xb_t(dx)\\
&\ascv&\phi(F)
\end{eqnarray*}
as $t\uparrow\eta$. Since there exists a countable separating set of
bounded Lipschitz functions $\{\phi_n\}$, and the above holds for
each $\phi_n$,
\[
F' = F+ \gamma\int_0^{\eta}\Zb_s \,ds\qquad \mbox{a.s. }
\]
\upqed\end{pf*}

\begin{remark}
(a) Theorem~\ref{TExtConv} holds in the critical branching
case. That is, if $\beta= 0$,
\[
X_t\ascv\delta_{F} \quad\mbox{and}\quad X_t'\ascv\delta_{F'},
\]
where
\[
F' =F+\int_0^\Tx\Zb_s\,ds.
\]
The convergence of the critical ordinary SOU process to a random point
follows directly from Tribe's result. That this holds for the SOU
process with attraction to the COM follows from the calculations above.
\begin{longlist}[(b)]
\item[(b)] The distribution of the random point $F$ has been identified
in Tribe~\cite{Tribe1992} by approximating with branching particle
systems. In fact, the law of $F$ can be identified as $x_\Tx$, where
$x_t$ is an Ornstein--Uhlenbeck process with initial distribution given
by $\Xb_0$ and $\Tx$ is the extinction time. Finding the distribution
of $F'$ remains an open problem however.
\end{longlist}
\end{remark}
%

\subsection{On the survival set, the attractive case}\label{SsSurv}
\mbox{}
\begin{pf*}{Proof of Theorem~\ref{TSurvival}}
Let $\gamma\ge0$ and as in the proof of Theorem~\ref{TExtinc}, assume
$\P(\tK_0(\vert y_0\vert^2))<\infty$. Also, without loss of
generality, assume that $d=1$ for this proof.

By Theorem~\ref{TIto}, $\Yb_t$ is a continuous local martingale with
decomposition given by $\Yb_t = \Yb_0 + M_t(Y)$ where
\[
[M(Y)]_t = \int^t_0\frac{\overline{Y^2_s}-\Yb_s^2}{K_s(1)}\,ds = \int
^t_0\frac{V(Y_s)}{K_s(1)}\,ds,
\]
with $V(Y_s) = \overline{Y^2_s}-\Yb_s^2$. Theorem IV.34.12 of \cite
{RW1985} shows that on the set $\{[M(Y)]_\infty< \infty\}\cap S$,
$M_t(Y)$ a.s. converges. 

Note that by Lemma~\ref{LExtProb}, for a.e. $\omega\in S$,
$W(\omega
)> 0$, recalling that $W= \lim_{t\rightarrow\infty} e^{-\beta
t}K_t(1)$. Hence, it follows that $[M(Y)]_\infty<\infty$ on $S$ if
\[
\int_0^\infty e^{-\beta s}V(Y_s)\,ds <\infty.
\]
%
Then
\begin{eqnarray*}
\P\biggl(\int_0^\infty e^{-\beta s}V(Y_s) \,ds; S \biggr) &\le&\P
\biggl(\int
_0^\infty e^{-\beta s} \overline{Y^2_s} \,ds ; S\biggr) \\
&\le&\int_0^\infty e^{-\beta s}A(\gamma, s) \,ds \\
&<& \infty
\end{eqnarray*}
by Cauchy--Schwarz and Lemma~\ref{Lsecmom} since $\gamma\ge0$.
Therefore on $S$, $\bar{Y}_t$ converges a.s. to some limit $\Yb
_\infty
$. Note that if $\gamma= 0$, Remark~\ref{RCOMrel}(b) gives $\Yb_t =
\Zb_t$ and so (b) holds.

That $\Zb_t$ converges on $S$ for $\gamma>0$ follows from the fact that
$\Yb_t$ converges and equation (\ref{Ecorr}) by setting
\begin{eqnarray*}
\Zb_t &=& \Yb_t - \gamma\int_0^t e^{-\gamma(t-s)}\Yb_s \,ds \\
&=& \Yb_t -\Yb_\infty+ \gamma\int_0^t e^{-\gamma(t-s)}(\Yb
_\infty-\Yb _s) \,ds +e^{-\gamma t}\Yb_\infty\\
&\rightarrow&0\qquad \mbox{as } t\rightarrow\infty.
\end{eqnarray*}

By Remark~\ref{RCOMrel}(b), we see that for $\gamma>0$ since $\Yb_t
= \Zb_t + \gamma\int_0^t \Zb_s \,ds$, $\Zb_t \ascv0$ and $ \Yb_t
\ascv
\gamma\int_0^\infty\Zb_s \,ds.$

Now argue by conditioning as in the end of Theorem~\ref{TExtinc} to
get the full result.
\end{pf*}

The next few results are necessary to establish the almost sure
convergence of~$\Xb_t$ on the survival set. This will in turn be used
to show the almost sure convergence of $\Xb_t'$ using the
correspondence of Remark~\ref{RCOMrel}(a).

Let $P_t$ be the standard Ornstein--Uhlenbeck semigroup (with attraction
to the origin). Note that $P_t \rightarrow P_\infty$ in norm where
\[
P_\infty\phi(x) = \int\phi(z)\biggl(\frac{\gamma}{\pi}\biggr)^{{d}/{2}}
e^{-\gamma\vert z\vert^2}\,dz,
\]
which is independent of $x$. Recall that $W = \lim_{t\rightarrow
\infty
} e^{-\beta t}X_t(1)$ and $S=\{W>0\}$ a.s. from Lemma~\ref{LExtProb}.

\begin{lemma}\label{LSOUconvL2}
If $\gamma>0$, $\P(X_0(\vert x\vert^4))<\infty$ and $\P
(X_0(1)^4)<\infty$, then on $S$, for any $\phi\in\Lipp$, $e^{-\beta
t}X_t(\phi)\ltwocv WP_\infty\phi$ and
\[
\P\bigl(\vert e^{-\beta t}X_t(\phi) - WP_\infty\phi\vert^2
\bigr)\le
Ce^{-\zeta t},
\]
where $C$ depends only on $d$ and $X_0$, and $\zeta$ is a positive
constant dependent only on $\beta$ and $\gamma$.
\end{lemma}


\begin{remark} \label{Rconv}
As the $L^2$ convergence in Lemma~\ref{LSOUconvL2} is exponentially
fast, it follows from the Borel--Cantelli lemma and Chebyshev inequality
that for a strictly increasing sequence $\{t_n\}_{n=0}^\infty$ where
$|\{t_n\}\cap[k, k+1)| = \lfloor e^{\zeta k/2}\rfloor$, for
$\phi\in\Lipp$,
\[
e^{-\beta t_n}X_{t_n}(\phi) \rightarrow WP_\infty\phi\qquad\mbox{a.s. as }
n\rightarrow\infty.
\]
\end{remark}

The idea is to use the above remark to bootstrap up to almost sure
convergence in Lemma~\ref{LSOUconvL2} with some estimates on the
modulus of continuity of the process $e^{-\beta t}X_t(\phi).$

\begin{lemma} \label{LSOU-increments}
Suppose $\gamma>0$, $\P(X_0(|x|^8))<\infty$ and $\P
(X_0(1)^8)<\infty$. If $\phi\in\Lipp$ and $h>0$, then
%
\begin{equation}\label{EIncrBd}
\P\bigl(\bigl[e^{-\beta(t+h)}X_{t+h}(\phi) - e^{-\beta t}X_t(\phi)
\bigr]^4\bigr) \le C(t)h^2e^{-\zeta^*t},
\end{equation}
where $\zeta^*$ is a positive constant depending only on $\beta$ and
$\gamma$ and $C$ is polynomial in~$t$, and depends on $\gamma$,
$\beta$
and $d$.
\end{lemma}


Let $\Psi\dvtx \R^d\rightarrow\R$, $ p\dvtx \R_+\rightarrow\R$ be positive,
continuous functions. Further, suppose $\Psi$ is symmetric about 0 and
convex with $\lim_{|x|\rightarrow\infty} \Psi(x) = \infty$ and $p(x)$
is increasing with $p(0) = 0$. The following is a very useful result of
Garsia, Rodemich and Rumsey~\cite{GRR1970}.

\begin{proposition} \label{PGRR}
If $f$ is a measurable function on $[0,1]$ such that
%
\begin{equation}\label{EGRR1}
\int\int_{[0,1]^2} \Psi\biggl(\frac{f(t)-f(s)}{p(|t-s|)}\biggr) \,ds \,dt = B
<\infty
\end{equation}
then there is a set $K$ of measure 0 such that if $s, t\in
[0,1]\setminus K$ then
%
\begin{equation}\label{EGRR2}
\vert f(t)-f(s)\vert \le8 \int^{|t-s|}_0 \Psi^{-1}\biggl(\frac{B}{u^2}\biggr)\,dp(u).
\end{equation}
If $f$ is also continuous, then $K$ is in fact the empty set.
\end{proposition}

With this result in hand, we can now bring everything together to prove
convergence of $\Xb_t$.

%
\begin{pf*}{Proof of Theorem \protect\ref{TSOUconv}}
The strategy for this proof is simple: We use Remark~\ref{Rconv} to
see that we can lay down an increasingly (exponentially) dense sequence
$e^{-\beta t_n}X_{t_n}(\phi)$ which converges almost surely, and that
we can use Lemma~\ref{LSOU-increments} to get a modulus of continuity
on the process $e^{-\beta t} X_{t}(\phi)$, which then implies that if
the sequence is converging, then the entire process must be converging.

Assume that $\P(K_0(|y_0|^8))=\P(X_0(|x|^8))<\infty$
and $\P
(X_0(1)^8)<\infty$ and argue as in Theorem~\ref{TExtinc} in the
general case. Let $\phi\in\Lipp$. Denote $ e^{-\beta t}X_t$ by
${\mathring{X}}
_t$ for the remainder of the proof. Let $T>0$, and let $\Psi(x) =
|x|^4$ and $p(t) = |t|^{3/4} (\log(\frac{\lambda}{t}))^{1/2}$
where $\lambda= e^4$. Let $B_T(\omega)$ be the constant\vspace*{-2pt} $B$ that
appears in Proposition~\ref{PGRR}, with aforementioned functions
$\Psi
$ and $p$, for the path ${\mathring{X}}_{Tt}(\omega), t\in[0,1]$.

Then note that
\begin{eqnarray*}
\P(B_T) &\equiv&\P\biggl[\int\int_{[0, 1]^2} \Psi\biggl(\frac {{\mathring
{X}}_{Tt}-{\mathring{X}} _{Ts}}{p(|t-s|)}\biggr) \,ds \,dt\biggr] \\
&=&\int\int_{[0, 1]^2} \frac{ \P[\vert{\mathring
{X}}_{Tt}-{\mathring{X}} _{Ts}\vert^4]}{|t-s|^3\log^2({\lambda
}/{|t-s|}) }\,ds \,dt \\
&\le&\int\int_{[0, 1]^2} \frac{ C(T(s\wedge t))e^{-\zeta
^*(s\wedge
t)}T^2|t-s|^2}{|t-s|^3\log^2({\lambda}/{|t-s|}) }\,ds \,dt \\
&=& 2T^2 \int_0^1\int_0^t \frac{ C(Ts)e^{-\zeta
^*(Ts)}|t-s|^2}{|t-s|^3\log^2({\lambda}/{|t-s|}) }\,ds \,dt \\
&\le&2C(T)T^2 \int_0^1\int_0^t \frac{ 1}{|t-s|\log^2(
{\lambda }/{|t-s|}) }\,ds \,dt \\
&\le&\frac{C(T)T^2}{2e^4},
\end{eqnarray*}
%
where $C$ is the polynomial term that appears in Lemma~\ref{LSOU-increments}. Since ${\mathring{X}}_t$ is continuous, by
Garsia--Rodemech--Rumsey~\cite{GRR1970}, for all $s, t \le1$,
\begin{eqnarray*}
\vert{\mathring{X}}_{Tt} - {\mathring{X}}_{Ts}\vert &\le&8 \int
_0^{\vert t-s\vert} \biggl(\frac {B_T}{u^2}\biggr)^{{1}/{4}}\,dp(u) \\
&\le& AB_T^{{1}/{4}}|t-s|^{{1}/{4}}\biggl(\log\frac{\lambda
}{|t-s|}\biggr)^{{1}/{2}},
\end{eqnarray*}
where $A$ is a constant independent of $T$ (see Corollary 1.2 of
Walsh~\cite{W1984} for this calculation). Rewriting the above,
%
\begin{equation}\label{EMC}
\vert{\mathring{X}}_t - {\mathring{X}}_s\vert \le D_T|t-s|^{{1}/{4}}\biggl(\log\frac{\lambda T}{|t-s|}\biggr)^{{1}/{2}}\qquad \forall s<t\le T,
\end{equation}
where $D_T \equiv A(\frac{B_T}{T})^{{1}/{4}}$. Note that $\P
(D_T^4) =\frac{ A^4TC(T)}{2e^4}$, which is\vspace*{2pt} polynomial in $T$ of fixed
degree $d_0>1$. Let $\Omega_0$ be the set of probability 1 such that
for all positive integers $T$ equation \eqref{EMC} holds and $D_T\le
T^{d_0} $ for $T$ large enough. To see \mbox{$P(\Omega_0)=1$}, use Borel--Cantelli:
\begin{eqnarray*}
\P(D_T\ge T^{d_0}) &=& \P(D_T^4\ge T^{4d_0})\\
&\le&\frac{\P(D_T^4)}{T^{4d_0}}\\
&\le&\frac{c}{T^{3d_0}}
\end{eqnarray*}
which is summable over all positive integers $T$.

Suppose $\omega\in\Omega_0$. Let $\delta^T(\omega)$ be such that
$\delta^{{-1}/{8}}(\log\frac{\lambda}{\delta})^{
{-1}/{2}} =
T^{d_0}$. Then for all integral $T>T_0(\omega)$, and $s, t\le T$ with
$|t-s|\le\delta$, $ \vert{\mathring{X}}_t - {\mathring{X}}_s\vert
\le
|t-s|^{{1}/{8}}.$

Now let $\{{\mathring{X}}_{t_n}\}$ be a sequence of the form in
Remark~\ref{Rconv}, with the additional condition that $\{t_n\}\cap[k, k+1)$ are
evenly spaced within $[k, k+1)$ for each $k\in\Z_+$ (i.e., $t_{n+1} -
t_n = ce^{-\zeta k/2}$ for $t_n\in\{t_n\}\cap[k, k+1)$). Evidently
${\mathring{X}}_{t_n}$ converges a.s. to a limit ${\mathring
{X}}_\infty$. Without loss of
generality, assume convergence of the sequence on the set $\Omega_0$.

There exists $T_1(\omega)$ such that for all $T>T_1$, $ce^{-\zeta T/2}
< \delta$. Hence, for all $t$ such that $T_1\vee T_0< t \le T$ there
exists $t'_n\in\{t_n\}$ such that $\vert t-t_n'\vert<ce^{-\zeta
\lfloor
t\rfloor/2}<\delta$ and hence
\begin{eqnarray*}
\vert{\mathring{X}}_t - {\mathring{X}}_\infty\vert &\le&\vert
{\mathring{X}}_t - {\mathring{X}}_{t_n'}\vert + \vert{\mathring
{X}} _{t_n'} - {\mathring{X}}_\infty\vert \\
&\le&\vert t-t_n'\vert^{{1}/{8}} + \vert{\mathring{X}}_{t_n'}-
{\mathring{X}}_\infty\vert \\
&\le& ce^{-{\zeta\lfloor t\rfloor}/{16}} + \vert{\mathring
{X}}_{t_n'} - {\mathring{X}} _\infty\vert.
\end{eqnarray*}
Sending $t\rightarrow\infty$ gives almost sure convergence of $
e^{-\beta t}X_t(\phi)$ to ${\mathring{X}}_\infty= WP_\infty(\phi)$
by Theorem \ref
{LSOUconvL2}, since $t_n' \rightarrow\infty$ with $t$. Note that
this implies for $\phi\in\Lipp$
%
\begin{equation} \label{Etestconv}
\Xb_t(\phi)\ascv P_\infty(\phi)
\end{equation}
since on $S$, $\frac{e^{\beta t}}{K_t(1)} \rightarrow W^{-1}$ a.s.

By Exercise 2.2 of~\cite{EK1986} on $\M_1(\R^d)$, the space of
probability measures on $\R^d$, the Prohorov metric of weak convergence
is equivalent to the Vasserstein metric.
It is easy to construct a class $\Theta$ that is a countable algebra of
Lipschitz functions and that therefore is strongly separating (see page~113 of~\cite{EK1986}).
Hence by Theorem~3.4.5(b) of~\cite{EK1986},
$\Theta$ is convergence determining. Since there exists a set
$S_0\subset S$ with $\P(S\setminus S_0) = 0$ such that on $S_0$,
equation \eqref{Etestconv} holds simultaneously for all $\phi\in
\Theta$,
\[
\Xb_t(\cdot)\rightarrow P_\infty(\cdot)
\]
in the Vasserstein metric, for $\omega\in S_0$ because $\Theta$ is
convergence determining.

To drop the dependence on the eighth moment, we argue as in the proof
of Theorem~\ref{TExtinc}, where we make use of the Markov Property and
the Compact support property for Historical Brownian Motion.
\end{pf*}


\begin{pf*}{Proof of Theorem~\ref{TISOUConv}}
This follows almost immediately from Theorem~\ref{TSOUconv} and the
representation given in Remark~\ref{RCOMrel}(a). Let $\phi\in\Lipp$, then
\begin{eqnarray*}
\Xb'_t(\phi) &=& \tK_t(\phi(Y_t)) =\tK_t\biggl(\phi\biggl(Z_t + \gamma \int
_0^t\Zb_s \,ds \biggr)\biggr)\\
&= &\int\phi\bigl(x + f(t,\omega) \bigr)\,d\Xb_t(dx),
\end{eqnarray*}
where $f(t)= \gamma\int_0^t \Zb_s \,ds$. Remark~\ref{RCOMrel}(b) gives
$f(t) = \Yb_t - \Zb_t$, and hence $f(t) \ascv\Yb_\infty$ follows from
Theorem~\ref{TSurvival}. Note that
\begin{eqnarray*}
&&\vert\Xb'_t(\phi) - P_\infty^{\Yb_\infty}(\phi)\vert\\
&&\qquad \le\biggl\vert
\int\phi\bigl(x + f(t) \bigr)\,d\Xb_t(dx) - \int\phi\bigl(x + f(\infty) \bigr)\,d\Xb
_t(dx) \biggr\vert \\
&&\qquad\quad{}+\biggl\vert\int\phi\bigl(x + f(\infty) \bigr)\,d\Xb_t(dx) - \int\phi\bigl(x +
f(\infty) \bigr)\,d\Xb_\infty(dx) \biggr\vert \\
&&\qquad\le\vert f(t) - f(\infty)\vert + d(\Xb_t, \Xb_\infty)
\end{eqnarray*}
since $\phi\in\Lipp$. Taking the supremum over $\phi$ and the previous
theorem give
\begin{eqnarray*}
d(\Xb'_t, \Xb'_\infty) \le\vert f(t, \omega) - f(\infty, \omega
)\vert +
d(\Xb_t, \Xb_\infty) \ascv0.
\end{eqnarray*}
\upqed\end{pf*}

\subsection{The repelling case, on the survival set}\label{SsRep}

Much less can be said for the SOU process repelling from its center of
mass than in the attractive case. We can, however, show that the center
of mass converges, provided the rate of repulsion is not too strong,
which we recall was the first step toward showing the a.s. convergence of
the normalized interacting SOU process in the attractive
case. The situation here is more complicated since we prove that the
COM of the ordinary SOU process with repulsion diverges almost surely,
implying that results for convergence of $\Xb'$ will not simply be
established through the correspondence. We finish with some conjectures
on the limiting measure for the repelling case.

As in the previous section, assume $Z_0=Y_0=y_0$, unless stated otherwise.

\begin{pf*}{Proof of Theorem~\ref{TrCOM}}
Assume that $\P(\Kb(|y_0|^2)) <\infty$, like in the proof of
Theorem~\ref{TExtinc}. As in that theorem, this condition can be
weakened to just the finite initial mass condition using similar reasoning.

For part (a), note that as in Theorem~\ref{TSurvival}, $\Yb_t$ will
converge if $\P([\Yb]_t) < \infty$ and which holds if the following
quantity is bounded:
\begin{eqnarray*}
\P\biggl(\int_0^\infty\frac{\overline{\vert Y_s\vert^2} - \Yb
_s^2}{e^{\beta
s}} \,ds; S \biggr) &\le&\P\biggl(\int_0^\infty\frac{\overline
{Y^2_s}}{e^{\beta s}} \,ds ; S\biggr) \\
&\le& c\int_0^\infty\frac{1+s^6e^{-2\gamma s}}{e^{\beta s}} \,ds \\
&<& \infty,
\end{eqnarray*}
by Lemma~\ref{Lsecmom} and by the conditions on $\gamma$.

For (b), we require the following lemma.

\begin{lemma}Let $-\beta/2<\gamma<0$ and $X_0\neq0$.
For a measure $m$ on $\R^d$, let $\tau_a(m)$ be $m$ translated by
$a\in
\R^d$. That is, $\tau_a(m)(\phi) = \int\phi(x+a)m(dx)$. Then:
\begin{longlist}[(ii)]
\item[(i)] For all but at most countably many $a$,
%
\begin{equation}\label{Enzlim}
P_{\tau_a(X_0)}(e^{\gamma t}\bar Z_t\to L\neq0|S)=1.
\end{equation}
\item[(ii)] For all but at most one value of $a$
\[
P_{\tau_a(X_0)}(e^{\gamma t}\bar Z_t\to L\neq0|S)>0.
\]
\end{longlist}
\end{lemma}

\begin{pf}
We first note that, by the correspondence \eqref{Ecorr}, we have that
\[
\Zb_t = \Yb_t - \gamma e^{-\gamma t}\int_0^te^{\gamma s}\Yb_s \,ds.
\]
Under our hypotheses $\Yb_t$ converges by Theorem~\ref{TrCOM}(a),
and hence
%
\begin{equation}\label{Erep1}
\lim_{t\to\infty} e^{\gamma t}\bar Z_t+\int_0^t\gamma e^{\gamma
s}\bar Y_s\,ds=0
\end{equation}
a.s. on $S$. Therefore, on $S$,
%
\begin{equation}\label{Erep2}
\lim_{t\to\infty} e^{\gamma t}\bar Z_t\qquad \mbox{exists a.s.}
\end{equation}


Note that one can build a solution of $\mbox{(SE)}^2$ with initial
conditions given by $\tau_a(X_0)$ by seeing that if $Y_t$ gives the
solution of~\ref{SE2}, then $Y_t+ a$ gives the solution
of $\mbox{(SE)}^2_{Y_0+a, K}$, and that the projection
\[
X'_t(\cdot) = \int\1(Y_t+a \in\cdot) K_t(dy)
\]
gives the appropriate interacting SOU process.

By \eqref{Erep1} and \eqref{Erep2},
\begin{eqnarray*}
\P_{X_0}(\{e^{\gamma t}\Zb_t\rightarrow L \neq0\}^c | S ) &=&
\P_{\tau_a(X_0)}\biggl(\lim_{t\rightarrow\infty} \int_0^t e^{\gamma
s}\Yb _s\,ds = 0 \Big|S\biggr) \\
&=& \P_{X_0}\biggl(\lim_{t\rightarrow\infty} \int_0^t e^{\gamma s} \frac
{K_s(a+Y_s)}{K_s(1)} \,ds = 0\Big|S\biggr) \\
&=& \P_{X_0}\biggl(-\frac{a}{\gamma}+ \lim_{t\rightarrow\infty} \int
_0^t e^{\gamma s}\Yb_s \,ds = 0\Big|S\biggr).
\end{eqnarray*}


The random variable $\int_0^\infty e^{\gamma s}\bar Y_s\,ds$ is finite
a.s. and so only a countable number of values $a$ exist with the
latter expression positive, implying the first result. The second
result also follows as well since the last expression in the above
display can be 1 for at most 1 value of $a$.
\end{pf}

To complete the proof of Theorem~\ref{TrCOM}(b), choose a value
$a\in
\R^d$ such that~\eqref{Enzlim} holds. By Theorem III.2.2 of \cite
{Perkins2002} and
the fact that $X_0P_s \ll\tau_a(X_0)P_t, $ for all $0<s\le t$, for
the OU semigroup $P_t$, we have that for all $0<s\le t$
%
\begin{equation} \label{Eabscon}
\P_{X_0}(X_{s+\cdot} \in\cdot) \ll\P_{\tau_a(X_0)}(X_{t+\cdot}
\in\cdot).
\end{equation}
By our choice of $a$,
\[
\P_{\tau_a(X_0)}\Bigl(\P_{X_1}\Bigl(\lim_{t\rightarrow\infty} e^{\gamma
t}\Zb_t = 0, S \Bigr) \Bigr) = 0,
\]
holds, and hence by \eqref{Eabscon} we have
\[
\P_{X_0}\Bigl(\lim_{t\rightarrow\infty} e^{\gamma t}\Zb_t = 0, S \Bigr)
= \P
_{X_0}\Bigl(\P_{X_1}\Bigl(\lim_{t\rightarrow\infty} e^{\gamma t}\Zb _t = 0,
S \Bigr) \Bigr)= 0.
\]
Recalling from \eqref{Erep2} that $\lim_{t\to\infty}e^{\gamma
t}\bar
Z_t$ exists a.s., we are done.
\end{pf*}

Note that for $0 > \gamma> -\frac{\beta}{2}$, this implies that even
if mass is repelled at rate $\gamma$, the COM of the interacting SOU
process still settles down in the long run. That is, driving $Y_t$ away
from $\Yb_t$ seems to have the effect of stabilizing it. One can think
of this as a situation where the mass is growing quickly enough that
the law of large numbers overcomes the repelling force.

More surprising is that the COM of the ordinary SOU process diverges
exponentially fast, even while the COM of the interacting one settles
down. This follows from the correspondence
\[
\Yb_t = \Zb_t + \gamma\int_0^t \Zb_s \,ds,
\]
and the cancellation that occurs in it due to the exponential rate of
$\Zb_t$.


The next lemma shows that Theorem 1 of Engl\"ander and Winter \cite
{EW2006} can be reformulated to yield a result for the SOU process
with repulsion at rate~$\gamma$ (where~$\gamma$ is taken to be a
negative parameter in our setting).

\begin{lemma} \label{Lreform}
On $S$, for the SOU process, $X,$ with repulsion rate $-\frac{\beta
}{d}< \gamma<0 $ and compactly supported initial measure $\mu$, and
any $\psi\in C_c^+(R^d)$
\[
e^{d\vert\gamma\vert t}\Xb_t(\psi) \pcv\xi\int_{\R^d} \psi(x) \,dx,
\]
where $\xi$ is a positive random variable on the set $S$.
\end{lemma}

\begin{pf}
Note that by Example 2 of Pinsky~\cite{Pinsky1996} it is shown that
the hypotheses of Theorem 1 of~\cite{EW2006} hold for the SOU process
with repulsion from the origin at rate $0< -\gamma<\frac{\beta}{d}$.
The theorem says that there is a function $\phi_c \in C_b^\infty(\R
^d)$ such that
%
\begin{equation}\label{EEW}
\frac{X_t(\psi)}{\E^{\mu}(X_t(\psi))} \pcv\frac{W\xi}{\mu
(\phi_c)},
\end{equation}
where $W$ is as in Lemma~\ref{LExtProb}.

Example 2 also shows that for $\psi\in C^+_c(\R^d)$,
\[
\lim_{t\rightarrow\infty} e^{-(\beta+\gamma d)t}\E^\mu(X_t(\psi
)) =
\mu(\phi_c)m(\psi),
\]
where $m$ is Lebesgue measure on $\R^d$. Hence, manipulating the
expression in \eqref{EEW} by using the previous equation and
Lemma \ref
{LExtProb} gives
\begin{eqnarray*}
e^{\vert\gamma\vert \,d t}\Xb_t(\psi) &\pcv&\frac{\xi W}{\mu(\phi_c)}
\lim
_{t\rightarrow\infty} \frac{e^{\vert\gamma\vert \,dt} \E^\mu
(X_t(\psi ))}{X_t(1)} \\
&=& \frac{\xi W}{\mu(\phi_c)} \lim_{t\rightarrow\infty} \frac
{e^{-(\beta+ \gamma d)t} \E^\mu(X_t(\psi))}{e^{-\beta t}X_t(1)}
\\
&=& \frac{\xi W}{\mu(\phi_c)} \frac{\mu(\phi_c)m(\psi)}{W}.\hspace*{80pt}\qed
\end{eqnarray*}
\noqed\end{pf}

This lemma indicates that on the survival set, when $\gamma<0$, one
cannot naively normalize $X_t$ by its mass since the probability
measures $\{\Xb_t\}$ are not tight. That is, a proportion of mass is
escaping to infinity and is not seen by compact sets. Note that the
lemma above implies that for $X_t$, the right normalizing factor is
$e^{(\beta+ \gamma d)t}$.

\begin{definition}
We say a measure-valued process ${\mathring{X}}$ goes locally extinct
if for any
finite initial measure ${\mathring{X}}_0$ and any bounded $A\in\B(\R
^d)$, there
is a $\P_{{\mathring{X}}_0}$-a.s. finite stopping time $\tau_A$ so
that ${\mathring{X}}_t(A)
= 0$ for all $t\ge\tau_A$ a.s.
\end{definition}

\begin{remark} \label{RLocExtinc}
Example 2 of Pinsky~\cite{Pinsky1996} also shows that for $\gamma\le
-\beta/d$ the SOU undergoes local extinction (all the mass escapes).
Hence for $\psi\in C_c(\R^d)$, there is no normalization where
$X_t(\psi)$ can be expected to converge to something nontrivial.
\end{remark}

From Remark~\ref{RIto2}, one can show that
\[
\Zb_t = \Zb_0 + N_t - \gamma\int_0^t \Zb_s \,ds,
\]
where $N$ is a martingale. Therefore, you can think of the COM of $X$,
the SOU process repelling from origin, as being given by an exponential
drift term plus fluctuations. The correspondence of Remark~\ref{RCOMrel}(b) implies then that
$\Yb_t = \Yb_0 + N_t,$
or in other words, the center of mass of the SOU process repelling from
its COM is given by simply the fluctuations.


We finish with some conjectures.


\begin{conjecture} \label{Crepulse}
On the survival set, if $X'_0$ is fixed and compactly supported, then
the following is conjectured to hold:
\begin{longlist}[(a)]
\item[(a)] If $0 < -\gamma< \frac{\beta}{d}$, then there exists
constant $\beta+ \gamma d \le\alpha< \beta$ so that for $\phi\in
C_c(\R^d)$,
\[
e^{-\alpha t}X'_t(\phi)\pcv\nu(\psi),
\]
where $\nu$ is a random measure depending on $\Yb_\infty$.
\item[(b)] If $\beta/d \le-\gamma$, then $X'_t$ undergoes local extinction.
\end{longlist}
\end{conjecture}

We expect that $\alpha<\beta$ simply because of the repulsion from
the COM built in to the model results in a proportion of mass being
lost to infinity. One would expect that the limiting measure $\nu$ is a
random multiple of Lebesgue measure as in the ordinary SOU process
case, due to the correspondence, but it is conceivable that it is some
other measure which has, for example, a dearth of mass near the
limiting COM.

As stated earlier, it is difficult to use Lemma~\ref{Lreform} to prove
this conjecture as the correspondence becomes much less useful\vadjust{\goodbreak} in the
repulsive case. The problem is that while the equation
\[
\int\phi(x) \,dX'_t(x) = \int\phi\biggl(x+ \gamma\int_0^t \Zb_s \,ds\biggr)
\,dX_t(x)
\]
still holds for $t$ finite, the time integral of $\Zb_s$ now diverges.



\section{Proofs of technical lemmas}\label{Stech}

We now prove the lemmas first stated in Section~\ref{Sexpre}.

\begin{pf*}{Proof of Lemma~\ref{LExtProb}}
We first note that $\brK_t = e^{-\beta t}K_t$ is a $(\Delta/2, 0,\break
e^{-\beta t})$-historical process. The martingale problem then shows
that $\brK_t(1)$ is a nonnegative $(\F_t)$-martingale and therefore
converges almost surely by the Martingale Convergence theorem to a
random variable $W$. It follows that $\{\Tx< \infty\} \subset\{W=0\}
$, since~0 is an absorbing state for $K_t(1)$.
Exercise II.5.3 in~\cite{Perkins2002} shows that
\[
\P(\Tx< \infty) = e^{-2\beta K_0(1)}.
\]
The same exercise also shows
\[
\P_{K_0}(\exp{(-\lambda\brK_t(1)})) = \exp\biggl(-\frac{2\lambda
\beta
\brK_0(1)}{2\beta+ \lambda(1-e^{-\beta t})}\biggr).
\]
Now sending $t \rightarrow\infty$ gives
\[
\P_{ K_0}(e^{-\lambda W}) = \exp\biggl(-\frac{2\beta\lambda\brK
_0(1)}{2\beta+\lambda}\biggr),
\]
and sending $\lambda\rightarrow\infty$ gives
\[
\P(W=0) = e^{-2\beta\brK_0(1)} = \P(\Tx<\infty)
\]
since $K_0 = \brK_0$. Therefore, $\{\Tx< \infty\} = \{W=0\}$ almost surely.
\end{pf*}


\begin{pf*}{Proof of Lemma~\ref{LCSP}}
We follow the proof of Theorem III.1.3(a) in Perkins~\cite
{Perkins2002}. First, note
that if $H$ is another supercritical historical Brownian motion
starting at time $\tau$ with initial measure $m$ under $\Q_{\tau, m}$
and $A$ a Borel subset of~$C$, then the process defined by
\[
H'_t(\cdot) = H_t(\cdot\cap\{y\dvtx  y_\tau\in A\})
\]
is also a supercritical historical Brownian motion starting at time
$\tau$ with initial measure $m'$ given by $m'(\cdot)=m(\cdot\cap
A)$ under $\Q_{\tau, m}$. Then using the extinction probabilities for
$H'$ (refer, e.g., to Exercise II.5.3 of~\cite{Perkins2002})
we have
%
\begin{equation}\label{EExt}
\Q_{\tau, m}\bigl(H_t(\{y\dvtx y_\tau\in A\})=0\ \forall t\ge s\bigr) =
\exp
{\biggl\{-\frac{2\beta m(A)}{1-e^{-\beta(s-\tau) }}\biggr\} }.
\end{equation}
Using the Markov property for $K$ at time $\frac{j}{2^n}$ and \eqref
{EExt} gives
%
\begin{eqnarray} \label{EBC}
&&\P\biggl[\exists t > \frac{j+1}{2^{n}} \mbox{ s.t. } K_t\biggl(\biggl\{y\dvtx \biggl\vert
y\biggl(\frac{j}{2^n}\biggr) - y\biggl(\frac{j-1}{2^n}\biggr)\biggr\vert>c\biggl(\frac
{j}{2^n}\biggr)h(2^{-n} ) \biggr\}\biggr) >0\biggr] \nonumber\\[-2pt]
&&\qquad=\P\bigl[1- \exp\bigl\{-\bigl( 2\beta K_{{j}/{2^n}} \bigl(\bigl\{y\dvtx \bigl\vert y(
{j}/{2^n}) - y\bigl({(j-1)}/{2^n}\bigr)\bigr\vert\nonumber\\[-2pt]
&&\hspace*{169pt}\qquad{}>c( {j}/{2^n})h(2^{-n} ) \bigr\}
\bigr)\bigr)\nonumber\\[-2pt]
&&\hspace*{197pt}\qquad{}\big/( 1 - e^{-\beta/{2^n}  }) \bigr\}\bigr]
\nonumber
\\[-8pt]
\\[-8pt]
\nonumber
&&\qquad\le\P\biggl[\frac{ 2\beta K_{{j}/{2^n}} (\{y\dvtx \vert y(
{j}/{2^n}) - y({(j-1)}/{2^n})\vert>c({j}/{2^n})h(2^{-n} ) \})}{
1 - e^{-\beta/{2^n} } } \biggr]
\\[-2pt]
&&\qquad\le\P\biggl[\frac{ 2\beta K_{{j}/{2^n}} (\{y\dvtx \vert y(
{j}/{2^n}) - y({(j-1)}/{2^n})\vert>c({j}/{2^n})h(2^{-n} ) \})}{
({\beta}/{2^n}) -( {\beta^2}/{2^{2n+1} }) } \biggr]\nonumber\\[-2pt]
&&\qquad\le\frac{2^{n+1}}{1 - ({\beta}/{2^{n+1} }) }\nonumber\\[-2pt]
&&\qquad\quad{}\times\P(K_{{j}/{2^n}}(1))\hat{\P}_{{j}/{2^n}}\biggl(\biggl\vert y\biggl(\frac{j}{2^n}\biggr) -
y\biggl(\frac{j-1}{2^n}\biggr)\biggr\vert>c\biggl(\frac{j}{2^n}\biggr)h(2^{-n} )\biggr),\nonumber
\end{eqnarray}
%
where in the last step we have simply rearranged the constants and
multiplied and divided by the mean mass at time $j2^{-n}$ and used the
definition of the Campbell measure (Definition~\ref{Dcamp}).

Since under the normalized mean measure, $y$ is a stopped Brownian
motion by Remark~\ref{RSOUHrep}(b), we use tail estimates to see that
the last quantity in \eqref{EBC} is
\begin{eqnarray*}
&\le&\frac{2^{n+1}}{1 - ({\beta}/{2^{n+1} })} \P(K_{
{j}/{2^n}}(1)) c_d n^{d/2 - 1}2^{-nc({j}/{2^n})^2/2} \\[-2pt]
&\le&2^{n+2}2^{{\beta j\ln2}/{2^n}}\P(X_0(1))c_d n^{d/2 -
1}2^{-n(({j}/{2^n})+c_0)/2}.
\end{eqnarray*}
Hence, summing over $j$ from $1$ to $n2^n$ gives
\begin{eqnarray*}
&&\P\biggl[\mbox{There exists } 1\le j\le n2^n \mbox{ s.t. } \exists
t >
\frac{j+1}{2^{n}} \mbox{ s.t. } \\[-2pt]
&&\quad K_t\biggl(\biggl\{y\dvtx \biggl\vert y\biggl(\frac{j}{2^n}\biggr) - y\biggl(\frac {j-1}{2^n}\biggr)\biggr\vert>c\biggl(\frac
{j}{2^n}\biggr)h(2^{-n} ) \biggr\}\biggr)>0 \biggr] \\[-2pt]
&&\qquad\le\P(X_0(1))c_d n^{d/2 - 1}2^{n+2- nc_0/2}\sum_{j=1}^{n2^n}
2^{{\beta j\ln2}/{2^n} - n({j}/{2^n})}\\[-2pt]
&&\qquad\le\P(X_0(1))c_d n^{d/2 - 1}2^{-2n+2}\sum_{j=1}^{n2^n}
2^{-
{j}/{2^n} }\\[-2pt]
&&\qquad\le\P(X_0(1))c_d n^{d/2 - 1}2^{-n+2},
\end{eqnarray*}
where we have used the fact that $c_0>6$. Hence, the sum over $n$ of
the above shows by Borel--Cantelli that there exists for almost sure
$\omega$, $N(\omega)$ such that for all $n>N$, for all $1\le j\le
n2^n$, for all $t\ge\frac{j+1}{2^n}$, for $K_t$-a.a. $y$,
\[
\biggl\vert y\biggl(\frac{j-1}{2^n}\biggr) - y\biggl(\frac{j}{2^n}\biggr)\biggr\vert < c\biggl(\frac
{j}{2^n}\biggr)h(2^{-n} ).
\]
Letting $\delta(\omega) = 2^{-N(\omega)}$, note that on the dyadics, by
above, we have that
\begin{eqnarray*}
\P(\delta< \lambda) &=& \P\biggl(N > -\frac{\ln\lambda}{\ln2}\biggr) \\
&=& \P\biggl[\exists n> -\frac{\ln\lambda}{\ln2}, \exists j\le n2^n,
\exists t>\frac{j+1}{2^n}, \mbox{ s.t. } \\
&&\quad K_t\biggl(\biggl\{y\dvtx \biggl\vert y\biggl(\frac{j}{2^n}\biggr) - y\biggl(\frac {j-1}{2^n}\biggr)\biggr\vert>c\biggl(\frac
{j}{2^n}\biggr)h(2^{-n} ) \biggr\}\biggr)>0\biggr]\\
&\le&\sum_{n= \lfloor-{\ln\lambda}/{\ln2} \rfloor} C'(d,
c_0)n^{d/2 - 1}2^{2n-nc_0/2} \\
&\le &C(d, c_0, \varepsilon)\lambda^{({c_0}/{2})-\varepsilon},
\end{eqnarray*}
where $\varepsilon$ can be chosen to be arbitrarily small (though the
constant $C$ will increase as it decreases). The rest of the proof
follows as in Theorem III.1.3(a) of~\cite{Perkins2002}, via an
argument similar to
Levy's proof for the modulus of continuity for Brownian motion.
\end{pf*}


\begin{pf*}{Proof of Lemma~\ref{Lsecmom}}
Assume that $Z_0 = Y_0$ and $t<\Tx$. Recall
\[
Z_t(\omega, y) \equiv e^{-\gamma t}Z_0+ \int_0^t e^{\gamma(s-t)} \,dy_s.
\]
Note that below ``$\lesssim$'' denotes less than up to multiplicative
constants independent of $t$ and $y$. Suppose that $y\in S(\delta,
c(t))$, where $S(\delta, c(t))$ is the same as in the previous lemma.
Then, as $Y_t = Z_t + \gamma\int_0^t \bar{Z}_s\,ds,$
\[
|Y_t|^2 \lesssim|Z_t|^2 + \gamma^2t \int_0^t\vert\Zb_s\vert^2\,ds
\lesssim|Z_t|^2 + \gamma^2t \int_0^t \overline{|Z_s|^2}\,ds
\]
by Cauchy--Schwarz and Jensen's inequality. Therefore, integrating with
respect to the normalized measure gives
%
\begin{equation}\label{ECOMbd}
\overline{\vert Y_t\vert^2} \le\overline{|Z_t|^2} + \gamma^2t \int_0^t
\overline{|Z_s|^2}\,ds
\end{equation}
and therefore we need only find the appropriate bounds for expectation
of $\overline{|Z_t|^2}$ to get the result.

After another few applications of Cauchy--Schwarz and integrating by parts,
\begin{eqnarray*}
|Z_t|^2&\lesssim& e^{-2\gamma t}|Z_0|^2 + \biggl|e^{-\gamma t} \int_0^t
e^{\gamma s}\,dy_s \biggr|^2\\
&\lesssim& e^{-2\gamma t}|Z_0|^2 + e^{-2\gamma t}\biggl|e^{\gamma t} y_t
- y_0 -\gamma\int_0^t y_s e^{\gamma s}\,ds \biggr|^2\\
&\lesssim& e^{-2\gamma t}(|Z_0|^2+|y_0|^2) +|y_t|^2 + \gamma^2
t\int
_0^t |y_s|^2 e^{-2\gamma(t-s)}\,ds.
\end{eqnarray*}
As $y\in S(\delta, c(t))$,
\begin{eqnarray*}
|Z_t|^2&\lesssim& e^{-2\gamma t}(|Z_0|^2+|y_0|^2) + |y_0|^2 +
\biggl(\frac{tc(t)h(\delta)}{\delta}\biggr)^2 \\
&&{}+ \gamma^2 t\int_0^t \biggl[|y_0|^2+ \biggl(\frac
{sc(s)h(\delta)}{\delta}\biggr)^2\biggr] e^{-2\gamma(t-s)}\,ds \\
&\lesssim& e^{-2\gamma t}(|Z_0|^2+|y_0|^2) + |y_0|^2 + |y_0|^2\gamma
t(1-e^{-2\gamma t})/2 \\
&&{}+ c(t)^2 \biggl(\frac{h(\delta)}{\delta}\biggr)^2\bigl(t^2
+ \gamma t^3(1-e^{-2\gamma t})\bigr) \\
& \lesssim&(1+\vert\gamma\vert t)(1+e^{-2\gamma t})(|Z_0|^2+|y_0|^2)
\\
&&{} +c(t)^2 \biggl(\frac{h(\delta)}{\delta}\biggr)^2\bigl(t^2
+ \gamma t^3(1-e^{-2\gamma t})\bigr).
\end{eqnarray*}
Integrating by the normalized measure $\tK_t$,
\begin{eqnarray*}
\overline{|Z_t|^2} &\lesssim&(1+\vert\gamma\vert t)(1+e^{-2\gamma
t})\tK_t(|Z_0|^2+|y_0|^2) \\
&&{}+ c(t)^2 \biggl(\frac{h(\delta)}{\delta}\biggr)^2\bigl(t^2 + \gamma
t^3(1-e^{-2\gamma t})\bigr).
\end{eqnarray*}
%
Then using \eqref{ECOMbd} and using the above bound on $\overline
{|Z_t|^2}$ gives
%
\begin{eqnarray}\label{El2bd}
\overline{\vert Y_t\vert^2} &\le&\overline{|Z_t|^2} + \gamma^2t
\int_0^t
\overline{|Z_s|^2}\,ds \nonumber\\
&\lesssim&(1+t^2)(1+e^{-2\gamma t})\biggl(\tK_t(|Z_0|^2+|y_0|^2)+ \int
_0^t\tK_s(|Z_0|^2+|y_0|^2)\,ds \biggr) \\
&&{}+ (1+t)c(t)^2\biggl(\frac{h(\delta)}{\delta
}
\biggr)^2\bigl(t^2 + \gamma t^3(1-e^{-2\gamma t})\bigr). \nonumber
\end{eqnarray}
Note that $\phi(y) \equiv|Z_0(y_0)|^2+|y_0|^2$ and $\phi_n(y)\equiv
\phi(y)1(\vert y\vert\le n)$ are $\hat{\F}_0$ measurable. By applying
It\^
o's formula to $K_t(\phi_n)K_t^{-1}(1)$ and using the decomposition
$K_t(\phi_n) = K_0(\phi_n) + \int_0^t\phi_n(y)\,dM(s, y) + \beta\int_0^t
K_s(\phi_n) \,ds$ (which follows from Proposition 2.7 of \cite
{Perkins1995}), we get
\[
\tK_t(\phi_n) = \tK_0(\phi_n) + N_t(\phi_n),
\]
where $N_t(\phi_n)$ is a local martingale until time $\Tx$, for each
$n$. In fact, the sequence of stopping times $\{T_N\}$ appearing in
Theorem~\ref{TIto} can be used to localize each $N_t(\phi_n)$.
Applying first the monotone convergence theorem and then localizing gives
\begin{eqnarray*}
\P\bigl(\tK_t(\phi);t<\Tx\bigr) &=& \lim_{n\rightarrow\infty} \P\bigl(\tK
_t(\phi _n);t<\Tx\bigr) \\
&=& \lim_{n\rightarrow\infty}\lim_{N\rightarrow\infty} \P\bigl(\tK
_t(\phi _n)\1(t<T_N) \bigr) \\
&=& \lim_{n\rightarrow\infty}\lim_{N\rightarrow\infty} \P\bigl(\tK
_{t\wedge T_N}(\phi_n) - \tK_{T_N}(\phi_n)\1(t\ge T_N) \bigr)\\
&\le&\lim_{n\rightarrow\infty}\lim_{N\rightarrow\infty} \P(\tK
_{t\wedge T_N}(\phi_n)) \\
&=& \lim_{n\rightarrow\infty}\lim_{N\rightarrow\infty} \P\bigl(\tK
_0(\phi _n) + N_{t\wedge T_N}(\phi_n) \bigr)\\
&=& \lim_{n\rightarrow\infty} \P(\tK_0(\phi_n)) \\
&=& \P(\tK_0(\phi)),
\end{eqnarray*}
where we have used the positivity of $\phi_n$ to get the fourth line
and the monotone convergence theorem in the last line.
Further, note that
\begin{eqnarray*}
\P\biggl(\int_0^t\tK_s(\phi)\,ds ;t<\Tx\biggr) &\le&\P\biggl(\int_0^t\tK _s(\phi)
\1 (s<\Tx) \,ds\biggr) \\
&=& \int_0^t \P\bigl(\tK_s(\phi);s<\Tx\bigr) \,ds \\
&\le& t \P(\tK_0(\phi)),
\end{eqnarray*}
by the calculation immediately above. Thus, taking expectations
in \eqref{El2bd} and plugging in $c(t) = \sqrt{c_0+t}$ gives
\begin{eqnarray*}
\P(\overline{\vert Y_t\vert^2}; t<\Tx)\lesssim(1+t^3)(1+
e^{-2\gamma t})\P(\tK_0(\phi)) + t^5(1+ te^{-2\gamma t})\P\biggl(\frac
{h(\delta)^2}{\delta^2}\biggr).
\end{eqnarray*}
Now let $c_0$ be chosen so that $\supp(K_t) \subset S(\delta,c(t))$ and
$p_{c_0}(\lambda) = C\lambda^\alpha$ for $\lambda\in[0,1]$.
Note that
\begin{eqnarray*}
\P\biggl(\frac{h(\delta)^2}{\delta^2}\biggr)&=&\P\biggl(\frac{\ln
^+(1/\delta)}{\delta}\biggr) = \int_0^\infty\biggl(\frac{\ln
^+(1/\lambda
)}{\lambda}\biggr)\,d\P(\delta<\lambda) \\
&\le&\int_0^1 \biggl(\frac{\ln^+(1/\lambda)}{\lambda}\biggr)\,d\P
(\delta
<\lambda) +\P(\delta>1) \\
&=& \lim_{\lambda\downarrow0} \biggl( \frac{\ln^+(1/\lambda
)}{\lambda}\P
(\delta<\lambda)\biggr) - \P(\delta<1) \\
&&{}- \int_0^1 \P(\delta<\lambda)\,d\biggl(\frac{\ln
^+(1/\lambda
)}{\lambda}\biggr) +\P(\delta>1)\\
&<& \infty,
\end{eqnarray*}
by choosing the constant $c_0$ so that $\alpha$ is large ($\alpha\ge2$
is enough).
\end{pf*}

Remark~\ref{Rmoments} follows from the above proof, after noting that
the exponent $\alpha$ in Lemma~\ref{LCSP} can be made arbitrarily
large by choosing a sufficiently large constant~$c_0$. Hence by
choosing $\alpha$ appropriately, we can show that
\[
\P\biggl[\biggl(\frac{h(\delta)}{\delta}\biggr)^k\biggr] <\infty,
\]
which can then be used to adapt the proof above.

Recall that $\Lipp= \{\psi\in C(\R^d)\dvtx  \forall x, y, \vert\psi (x)
- \psi(y)\vert\le\vert x-y\vert, \|\psi\| \le1\}$. 
We will, with a slight abuse of notation, allow $M$ to denote the
orthogonal martingale measure generated by the martingale problem for
$X$. Let $A$ be the infinitesimal generator for an OU process, and
hence recall that for $\phi\in C^2(\R^d)$,
\[
A\phi(x) = -\gamma x\cdot\nabla\phi(x) +\frac{\Delta}{2}\phi(x).
\]

The next two proofs are for lemmas stated in Section~\ref{SsSurv}.

\begin{pf*}{Proof of Lemma~\ref{LSOUconvL2}}
Let $\phi\in\Lipp$.
By the extension of the martingale problem for $X$ given in Proposition
II.5.7 of~\cite{Perkins2002}, for functions $\psi\dvtx [0,T]\times\R
^d\rightarrow\R$ such
that $\psi$ satisfies the definition before that proposition,
\[
X_t(\psi_t) = X_0(\psi_0) +\int_0^t\int\psi_s(x)\,dM(x,s) + \int_0^t
X_s(A\psi_s + \beta\psi_s + \dot{\psi}_s) \,ds,
\]
where $M$ is the orthogonal martingale measure derived from the
martingale problem for the SOU process. It is not difficult to show
that $\psi_s = P_{t-s}\phi$ where $\phi$ as above satisfies
requirements for Proposition II.5.7 of~\cite{Perkins2002}. Plugging
this in gives
\[
X_t(\phi) = X_0(P_t\phi) + \int_0^t\int P_{t-s}\phi(x)\,dM(s, x) +
\int
_0^t\beta X_s(P_{t-s}\phi)\,ds
\]
since $\frac{\partial}{\partial s}P_s\phi= AP_s\phi$. Multiplying by
$e^{-\beta t}$ and integrating by parts gives
%
\begin{eqnarray} \label{EXblim}
e^{-\beta t}X_t(\phi) &=& e^{-\beta t}X_t(\psi_t)\nonumber\\
& =& X_0(\psi_0) -
\int
_0^t\beta e^{-\beta s}X_s(\psi_s)\,ds + \int_0^te^{-\beta s} \,dX_s(\psi
_s)\\
&= &X_0(P_t\phi) + \int_0^t\int e^{-\beta s}P_{t-s}\phi(x)\,dM(s,x).\nonumber
\end{eqnarray}
Note that as the $OU$-process has a stationary distribution $P_\infty$
where $P_t \rightarrow P_\infty$ in norm. When $s$ is large in \eqref
{EXblim}, $P_{t-s}\phi(x)$ does not contribute much to the stochastic
integral and hence we expect the limit of $e^{-\beta t}X_t(\phi)$ to be
%
\begin{equation}\label{EXblim2}
X_0(P_\infty\phi) + \int_0^\infty\int e^{-\beta s}P_\infty\phi
(x) \,dM(s,x),
\end{equation}
which is a well defined, finite random variable as 
%
\[
\biggl[\int_0^\cdot\int e^{-\beta s}P_\infty\phi(x) \,dM(s, x)\biggr]_\infty
< \|
\phi\|^2 \int_0^\infty e^{-2\beta s}X_s(1)\,ds,
\]
which is finite in expectation. As $P_\infty\phi(x)$ does not depend on
$x$, it follows that
\begin{eqnarray*}
\eqref{EXblim2} = (P_\infty\phi) X_0(1) + (P_\infty\phi)
\int
_0^\infty\int e^{-\beta s } \,dM(s, x) = WP_\infty\phi.
\end{eqnarray*}
Given this decomposition for $WP_\infty\phi$, we write
\begin{eqnarray*}
&&\P\bigl(\bigl(e^{-\beta t}X_t(\phi) - WP_\infty\phi\bigr)^2\bigr)\\
&&\qquad\le3 \P\biggl(\biggl(\int
_t^\infty e^{-\beta s}P_{\infty}\phi(x)\,dM(s,x) \biggr)^2 \biggr)\\
&&\qquad\quad{}+ 3\P\biggl(\biggl(\int_0^t\int e^{-\beta s}\bigl(P_{t-s}\phi(x)- P_{\infty }\phi
(x)\bigr)\,dM(s,x)\biggr)^2\\
&&\hspace*{160pt}\qquad{}+ X_0(P_\infty\phi- P_t\phi)^2\biggr).
\end{eqnarray*}

If $z_t$ is a $d$-dimensional OU process satisfying $dz_t = -\gamma z_t
dt+dB_t,$ where $B_t$ is a $d$-dimensional Brownian motion, then
\[
z_t =e^{-\gamma t}z_0+ \int_0^t e^{-(t-s)\gamma}\,dB_t
\]
and hence $z_t$ is Gaussian,\vspace*{-2pt} with mean $e^{-\gamma t}z_0$ and
covariance matrix $\frac{1}{2\gamma}(1-e^{-2\gamma t})I$. Evidently,
$z_\infty$ is also Gaussian, mean 0 and variance $\frac{1}{2\gamma} I$.
We use a simple coupling: suppose that $w_t$ is a random variable
independent of $z_t$ such that $z_\infty= z_t+ w_t$ (i.e., $w_t$ is
Gaussian with mean $-e^{-\gamma t}z_0$ and covariance $\frac
{1}{2\gamma
}e^{-2\gamma t}I$). Then using the fact that $\phi\in\Lipp$ and the
Cauchy--Schwarz inequality, followed by our coupling with $z_0 = x$ gives
\begin{eqnarray*}
X_0(P_\infty\phi- P_t\phi)^2 &=& \biggl( \int\E^x\bigl(\phi (z_\infty)
- \phi(z_t)\bigr)X_0(dx)\biggr)^2 \\
&\le&\int\E^x(|z_\infty- z_t|)^2X_0(dx)X_0(1) \\
&=& \int\E^x(|w_t|^2)X_0(dx)X_0(1) \\
&= &\int e^{-2\gamma t}\biggl(|x|^2 + \frac{d}{2\gamma}
\biggr)X_0(dx)X_0(1)\\
&\le& ce^{-2\gamma t}\biggl(\int\vert x\vert^2X_0(dx)X_0(1) + X_0(1)^2\biggr).
\end{eqnarray*}
%
Taking expectations and using Cauchy--Schwarz and the assumptions on
$X_0$ gives exponential rate of convergence for the above term.

Since we can think of $\int_0^r\int e^{-\beta s}P_{t-s}\phi(x) \,dM(s,x)$
as a martingale in $r$ up until time $t$, various martingale
inequalities can be applied to get bounds for the terminal element,
$\int_0^t\int e^{-\beta s}P_{t-s}\phi(x) \,dM(s,x)$. Note that this
process is not in general a martingale in $t$. Therefore, we have
%
\begin{eqnarray}\label{EM1}
&&\P\biggl[\biggl(\int_0^t\int e^{-\beta s}P_\infty\phi(x) \,dM(s,x) -
\int
_0^t\int e^{-\beta s}P_{t-s}\phi(x) \,dM(s,x) \biggr)^2 \biggr] \nonumber\\
&&\qquad= \P\biggl[\biggl(\int_0^t\int e^{-\beta s}\bigl(P_\infty
\phi
(x) - P_{t-s}\phi(x)\bigr) \,dM(s,x) \biggr)^2 \biggr] \\
&&\qquad\le\P\biggl[\int_0^t e^{-2\beta s}\int\bigl(P_\infty
\phi (x) - P_{t-s}\phi(x)\bigr)^2X_s(dx)\,ds\biggr]. \nonumber
\end{eqnarray}
Then as $\phi$ Lipschitz, by the coupling above,
\begin{eqnarray*}
\eqref{EM1} &\le&\P\biggl[\int_0^t e^{-2\beta s}\int e^{-2\gamma
(t-s)}\biggl(|x|^2 + \frac{d}{2\gamma} \biggr)X_s(dx)\,ds\biggr] \\
&= &\int_0^t e^{-2\beta s-2\gamma(t-s)}\P\biggl[\int|x|^2 + \frac
{d}{2\gamma} X_s(dx)\biggr]\,ds \\
&=& \int_0^t e^{-2\beta s-2\gamma(t-s)}\P\biggl[K_s(|Z_s|^2) +
\frac
{d}{2\gamma}X_s(1)\biggr]\,ds.
\end{eqnarray*}
Applying the Cauchy--Schwarz inequality followed by Remark \ref
{Rmoments}(b) gives
\begin{eqnarray*}
&&\P(K_s(|Z_s|^2))\\
&&\qquad=\P\bigl(\overline{|Z_s|^2}K_s(1) ; s<\Tx\bigr)
\\
&&\qquad\le\P(\overline{|Z_s|^2}^2; s<\Tx)^{{1}/{2}}\P
(X_s(1)^2)^{{1}/{2}}\\
&&\qquad\le B(s, \gamma, 4)^{{1}/{2}}\P(X_s(1)^2)^{{1}/{2}}\\
&&\qquad\le cB(s, \gamma, 4)^{{1}/{2}} e^{\beta s} \P\biggl(X_0(1)^2 + \frac
{1}{\beta}X_0(1)\biggr)^{{1}/{2}},
\end{eqnarray*}
where the last line follows by first noting that
\[
e^{-\beta t}X_t(1) = X_0(1) + \int_0^t\int e^{-\beta s} \,dM(s,x)
\]
is a martingale. That is,
%
\begin{eqnarray}\label{EmassL2}
e^{-2\beta s} \P(X_s(1)^2) &\le&2\P\biggl(X_0(1)^2 + \biggl(\int_0^s e^{-\beta
r} \,dM(r, x)\biggr)^2\biggr) \nonumber\\[-2pt]
&=& 2\P\biggl(X_0(1)^2+\biggl[\int_0^\cdot e^{-\beta r}\,dM(r, x)\biggr]_s\biggr)\nonumber
\\[-2pt]
&=& 2\P\biggl(X_0(1)^2+ \int_0^s e^{-2\beta r}X_r(1)\,dr \biggr)
\nonumber
\\[-8pt]
\\[-8pt]
\nonumber
&=& 2\P(X_0(1)^2)+ 2\int_0^s e^{-\beta r} \P(e^{-\beta r}X_r(1)
)\,dr\\[-2pt]
&=&2\P(X_0(1)^2)+ 2\int_0^s e^{-\beta r} \P(X_0(1) )\,dr\nonumber
\nonumber\\[-2pt]
&\le&2\P\biggl(X_0(1)^2 + \frac{1}{\beta}X_0(1)\biggr).\nonumber
\end{eqnarray}

Therefore,
\begin{eqnarray*}
\eqref{EM1} &\le&\int_0^t e^{-2\beta s-2\gamma(t-s)}
\biggl[e^{\beta
s}B(s, \gamma, 4)^{{1}/{2}} \P\biggl(X_0(1)^2 + \frac{1}{\beta
}X_0(1)\biggr)^{{1}/{2}}\biggr] \,ds\\[-2pt]
&&{}+ \int_0^te^{-2\beta s-2\gamma(t-s)} e^{\beta s}\frac
{d}{2\gamma}\P(X_0(1))\,ds \\[-2pt]
&\le&\int_0^t e^{-\beta s-2\gamma(t-s)}\biggl[B(s, \gamma, 4)^{
{1}/{2}} \P\biggl(X_0(1)^2 + \frac{1}{\beta}X_0(1)\biggr)^{{1}/{2}}\\[-2pt]
&&\hspace*{155pt}\qquad{} +
\frac
{d}{2\gamma}\P(X_0(1)) \biggr]\,ds \\[-2pt]
&<& C e^{-\zeta_1 t},
\end{eqnarray*}
where $\zeta_1 = \min(\beta, 2\gamma)-\varepsilon$ where
$\varepsilon$ is
arbitrary small and comes from the polynomial term in the integral.

Finally,
\begin{eqnarray*}
\P\biggl(\biggl(\int_t^\infty e^{-\beta s}P_{\infty}\phi(x)\,dM(s,x) \biggr)^2 \biggr)
&=&
(P_\infty\phi)^2\P\biggl(\int_t^\infty e^{-2\beta s}X_s(1) \,ds \biggr) \\[-2pt]
&\le&(P_\infty\phi)^2\int_t^\infty e^{-\beta s}\P(e^{-\beta
s}X_s(1))\,ds\\
&\le&\frac{(P_\infty\phi)^2}{\beta}e^{-\beta t}\P(X_0(1)),
\end{eqnarray*}
since $e^{-\beta s} X_s(1)$ is a martingale. Therefore, since $\zeta
_1<\beta$, we see that $\zeta= \zeta_1$ gives the correct exponent.
\end{pf*}

\begin{pf*}{Proof of Lemma~\ref{LSOU-increments}}
The proof will follow in a manner very similar to the proof of the
previous lemma. From the calculations above, we see that
\begin{eqnarray*}
&& e^{-\beta(t+h)}X_{t+h}(\phi) - e^{-\beta t}X_t(\phi)\\
&&\qquad =
X_0(P_{t+h}\phi- P_t\phi) \\
&&\qquad\quad{}+ \int_0^{t+h}\int e^{-\beta s}P_{t+h-s}\phi(x)\,dM(s, x) - \int
_0^{t}\int e^{-\beta s}P_{t-s}\phi(x)\,dM(s, x) \\
&&\qquad= X_0(P_{t+h}\phi- P_t\phi) + \int_0^{t}\int e^{-\beta s}
\bigl((P_{t+h-s}-P_{t-s})\phi(x) \bigr)\,dM(s, x)\\
&&\qquad\quad{}+ \int_t^{t+h}\int e^{-\beta s}P_{t+h-s}\phi(x)\,dM(s, x) \\
&&\qquad\equiv I_1+I_2 +I_3.
\end{eqnarray*}
Using the Cauchy--Schwarz inequality, we can find bounds for $\P
(|I_k|^4), k=1,2,3$, separately:
\begin{eqnarray*}
|I_1|^4 &=& X_0(P_{t+h}\phi- P_t\phi)^4 \\
&\le&\bigl[X_0\bigl((P_{t+h}\phi- P_t\phi)^2\bigr)X_0(1)\bigr]^2.
\end{eqnarray*}
Recalling the simple coupling in the previous lemma to see that
\begin{eqnarray*}
\bigl(P_{t+h}\phi(x) - P_t\phi(x)\bigr)^2 &\le&\E^x\bigl(|\phi(z_{t+h}) - \phi
(z_t)|^2\bigr) \\
&\le&\E^x(|z_{t+h} - z_t|^2)\\
&\le&\E^x(|w_{t,t+h}|^2),
\end{eqnarray*}
where $z$ is as above, an OU process started at $x$, and $w_{s,t}$ is
independent of $z_s$ but such that $z_t = z_s + w_{s,t}$. Hence,
$w_{s,t}$ is Gaussian with mean $x(e^{-\gamma t}-e^{-\gamma s})$ and
covariance matrix $\frac{I}{2\gamma}(e^{-2s\gamma}-e^{-2t\gamma})$.
Therefore,
\begin{eqnarray*}
\bigl(P_{t+h}\phi(x) - P_t\phi(x)\bigr)^2 &\le&|x|^2\bigl(e^{-\gamma
(t+h)}-e^{-\gamma t}\bigr)^2+ \frac{d}{2\gamma}\bigl(e^{-2t\gamma
}-e^{-2(t+h)\gamma}\bigr) \\
&=& e^{-2\gamma t}\biggl(|x|^2(1-e^{-\gamma h})^2+ \frac{d}{2\gamma
}(1-e^{-2h\gamma}) \biggr).
\end{eqnarray*}
Hence,
\begin{eqnarray*}
\P(|I_1|^4)
&=& \P\biggl[ \biggl(e^{-2\gamma t}(1-e^{-\gamma h})^2\int
|x|^2X_0(dx)X_0(1)\\
&&\hspace*{56pt}{} + \frac{d}{2\gamma}e^{-2\gamma t}(1-e^{2\gamma
h})X_0(1)^2\biggr)^2 \biggr]\\
&\le& C_1(d, \gamma) h^2 e^{-4\gamma t},
\end{eqnarray*}
%
where $C_1$ is a constant that is finite by assumptions on the initial measure.

To get bounds on the expectation of $I_2$, we use martingale
inequalities. Note that $ \int_0^{\cdot}\int e^{-\beta s}
((P_{t+h-s}-P_{t-s})\phi(x) )\,dM(s, x) =N(\cdot)$ is a martingale
until time $t$. Therefore, using the Burkholder--Davis--Gundy inequality
and the coupling above gives
\begin{eqnarray*}
\P(|I_2|^4) &\le& c\P([N]_t^2) \\
&=&c\P\biggl[\biggl(\int_0^t\int e^{-2\beta s}
\bigl((P_{t+h-s}-P_{t-s})\phi
(x) \bigr)^2X_s(dx) \,ds \biggr)^2\biggr]\\
&\le& c\P\biggl[\biggl( \int_0^t e^{-2\beta s} \bigl(e^{-\gamma
(t+h-s)}-e^{-\gamma(t-s)}\bigr)^2\int|x|^2X_s(dx)\,ds\biggr)^2 \\*
&&{}+ \biggl(\int_0^t\frac{de^{-2\beta s} }{2\gamma
}\bigl(e^{-2(t-s)\gamma}-e^{-2(t+h-s)\gamma}\bigr)X_s(1) \,ds \biggr)^2\biggr]\\
&\le& c\P\biggl[\biggl( t\int_0^t e^{-4\beta s} \bigl(e^{-\gamma
(t+h-s)}-e^{-\gamma(t-s)}\bigr)^4\biggl(\int|x|^2X_s(dx)\biggr)^2\,ds\biggr) \\*
&&{}+ t\biggl(\frac{d}{2\gamma}\biggr)^2\int_0^t
e^{-4\beta
s} \bigl(e^{-2(t-s)\gamma}-e^{-2(t+h-s)\gamma}\bigr)^2X_s(1)^2 \,ds \biggr]\\
&=&c t\int_0^t e^{-2\beta s} \bigl(e^{-\gamma(t+h-s)}-e^{-\gamma
(t-s)}\bigr)^4\P\biggl[\biggl(\int|x|^2e^{-\beta s}X_s(dx)\biggr)^2\biggr]\,ds \\
&&{}+ t\biggl(\frac{d}{2\gamma}\biggr)^2\int
_0^te^{-2\beta
s}\bigl(e^{-2(t-s)\gamma}-e^{-2(t+h-s)\gamma}\bigr)^2\P[(e^{-\beta
s}X_s(1))^2] \,ds.
\end{eqnarray*}
%
Since $X_s(|x|^2) = K_s( |Z_s|^2)$, by Remark~\ref{Rmoments}(b),
\begin{eqnarray*}
\P(e^{-2\beta s}X_s(|x|^2)^2) &\le&\P[ \overline
{Z^2_{s}}^2(e^{-\beta s}X_s(1))^2; s<\Tx]\\
&\le&\P( \overline{Z_{s}^2}^4; s<\Tx)^{1/2}\P
(e^{-4\beta s}X_{s}(1)^4)^{1/2} \\
&\le&\P( \overline{Z_{s}^8} ; s<\Tx)^{1/2} \P
(e^{-4\beta s}X_{s}(1)^4)^{1/2}\\
&\le& cB(\gamma, s, 8)^{1/2}\bigl( \P\bigl(X_0(1)^4+ sX_0(1)^2+ sX_0(1)\bigr)
\bigr)^{1/2}.
\end{eqnarray*}
The bound on the expectation of $e^{-4\beta s}X_s(1)^4$ follows by an
application of the BDG Inequality to $e^{-\beta s}X_s(1) =X_0(1) + \int
_0^s e^{-\beta r}\,dM(r, x)$ and similar calculations used to determine
the bound on \eqref{EmassL2}.
Therefore,
\begin{eqnarray*}
\P(|I_2|^4) &\le& ct\int_0^te^{-2\beta s}B(\gamma, s,
8)^{1/2} \bigl(e^{-\gamma(t+h-s)}-e^{-\gamma(t-s)}\bigr)^4\\
&&\hspace*{23pt}{}\times\bigl(\P\bigl(X_0(1)^4+ sX_0(1)^2 + sX_0(1)\bigr)
\bigr)^{1/2}\,ds \\
&&{}+ \frac{ct}{4\gamma^2} \int_0^t e^{-2\beta s}\bigl(e^{-2(t-s)\gamma
}-e^{-2(t+h-s)\gamma}\bigr)^2\P\biggl(X_0(1)^2 + \frac{1}{\beta}X_0(1)\biggr)\,ds
\\
&\le& ctB(\gamma, t, 8)^{1/2}\bigl(\P\bigl(X_0(1)^4+ tX_0(1)^2 +
tX_0(1)\bigr)\bigr)^{1/2}\\
&&{}\times (e^{-\gamma h}-1)^4\int_0^t e^{-2\beta
s}e^{-4\gamma(t-s)}\,ds \\
&&{}+ \frac{ct}{4\gamma^2}(e^{-2\gamma h} -1)^2\P\biggl(X_0(1)^2 + \frac
{1}{\beta}X_0(1)\biggr) \int_0^te^{-2\beta s}e^{-4(t-s)\gamma}\,ds \\
&\le& C_2(t,\gamma, \beta)h^2e^{-\zeta_1 t},
\end{eqnarray*}
where $C_2$ is polynomial in $t$. 
By another application of the BDG inequality, and noting that $\|\phi\|
= 1$,
\begin{eqnarray*}
\P(| I_3|^4) &=& \P\biggl[\biggl(\int_t^{t+h}\int e^{-\beta s} P_{t+h-s}\phi
(x)\,dM(s, x)\biggr)^4\biggr]\\
&\le& c\P\biggl[\biggl(\int_t^{t+h}\int\bigl(e^{-\beta s}P_{t+h-s}\phi (x)\bigr)^2X_s(dx)
\,ds\biggr)^2 \biggr] \\
&\le& c h \P\biggl[\int_t^{t+h}e^{-4\beta s}X_s(1)^2 \,ds \biggr]\\
&=& c he^{-2\beta t} \int_t^{t+h} \P[e^{-2\beta s}X_s(1)^2 ]\,ds \\
&\le& c he^{-2\beta t} \int_t^{t+h} \P[X_0(1)^2+X_0(1)/\beta]\,ds \\
&\le& C_3(\beta) h^2e^{-2\beta t},
\end{eqnarray*}
where the second last line follows from the same calculations performed
in estimating moments of $I_2$. Note that the constant $C_3$ does not
depend on $t$ here.

Putting the pieces together shows that there exists a function $C$
polynomial in $t$ and a positive constant $\zeta^*$ such that \eqref
{EIncrBd} holds.
\end{pf*}

\section*{Acknowledgement}
I am greatly indebted to my advisor, Ed Perkins, for his numerous (very
helpful) suggestions and generous support during the course of this
work. I am particularly grateful for his help with Part (b) of
Theorem~\ref{TrCOM}.

%


\printaddresses

\end{document}